# Asymptotic Efficiency of the Global Maximum-Likelihood Estimator in Multi-Emitter Localization Microscopy

**YI SUN***

*Electrical Engineering Department*
*The City College of City University of New York*
*New York, NY 10031, USA*
*Email: ysun@ccny.cuny.edu*
*ORCID: 0000-0002-2527-8311*

**Abstract:** In single-molecule localization microscopy (SMLM), a super-resolution image is built by estimating the positions of densely packed emitters from many noisy frames. We study the global maximum-likelihood (GML) estimator of the emitter positions for the multi-emitter SMLM imaging model with Poisson signal and mixed Poisson–Gaussian noise. We prove that, for any fixed configuration of resolvable emitters in 2D or 3D, the GML estimator is consistent, asymptotically normal, asymptotically unbiased, and asymptotically efficient, attaining the Cramér–Rao bound (CRB) as the number of frames grows in Theorem 1 or the average emitter intensity grows in Theorem 2. A corollary shows that in this limit the GML estimator and the unbiased Gaussian information-achieving estimator for a frame (UGIA-F) share the same limiting distribution, while their finite-information gap is the high-density advantage of a biased estimator. We realize the GML estimator by an expectation-maximization algorithm initialized near the true positions (EM-GML) and verify the theory by simulation, Theorem 1 with a 2D Gaussian point spread function and Theorem 2 with a 3D astigmatic one. In the finite-photon regime the biased EM-GML attains an error below the unbiased CRB, which motivates EM-GML as a benchmark and the design of learned localization algorithms for dense emitters.

## 1. Introduction

Single-molecule localization microscopy (SMLM) achieves spatial resolution far below the diffraction limit by recording many frames in which only a sparse, photoswitchable subset of emitters is active, estimating the position of each active emitter, and accumulating the positions into a super-resolution image [1-4]. The temporal resolution improves as more emitters are active per frame, but a higher density makes the point spread functions (PSFs) overlap, so the central estimation problem is the joint localization of several closely spaced emitters from a single noisy frame.

The accuracy of localization is limited by the Fisher information of the frame and the Cramér–Rao bound (CRB) it induces. For the SMLM imaging model with Poisson signal and mixed Poisson–Gaussian noise, Ref. [5] derived the closed-form Fisher information matrix for multiple emitters with arbitrary positions and the resulting CRB on the position estimates, Ref. [6] extended the analysis to unknown emitter intensities, and Ref. [7] introduced the root-mean-square minimum distance (RMSMD) metric together with the unbiased Gaussian information-achieving estimator for a frame (UGIA-F), which attains the per-frame Fisher information and serves as an unbiased benchmark [8-11].

The asymptotic behavior of the maximum-likelihood (ML) estimator is classical. The global maximum-likelihood (GML) estimator is consistent [12], asymptotically normal, and locally asymptotically efficient [13-15], as presented in modern form by Refs. [16, 17]. Its application to the multi-emitter SMLM model has not been made precise. Doing so would be significant since it would place the widely used ML localization on a rigorous estimation-theoretic footing, identify the exact accuracy limit, and establish when that limit is attained.

This paper makes that application precise. We define the GML estimator of the emitter positions for the SMLM Poisson working model and prove its asymptotic properties by verifying the regularity conditions of the model, in particular the strict positivity of the Fisher information for any fixed configuration of resolvable emitters, and then invoking the classical theory. Theorem 1 states that the GML estimator is consistent, asymptotically normal with covariance $(N\mathbf{F})^{-1}$, asymptotically unbiased, and asymptotically efficient as the number of frames $N$ grows, attaining the CRB where $\mathbf{F}$ is the Fisher information matrix of one frame. Theorem 2 states that for a single frame the GML estimator attains $\mathbf{F}(I)^{-1}$ as the average emitter intensity $I$ grows, again the CRB. A corollary shows that in either limit the GML estimator and the UGIA-F benchmark converge to the same distribution, so their finite-information gap is the high-density advantage studied as Hypothesis 1. The conditions are stated for a general PSF, so the results apply to the 2D Gaussian PSF [18], the 3D astigmatic PSF [19], the double-helix PSF [20], and other engineered or experimentally measured PSFs [21, 22].

We verify the theory with the expectation-maximization (EM) algorithm [23]. The EM algorithm decomposes the joint $M$-emitter localization into $M$ noise-free single-emitter ML problems, each solved by Newton iteration or gradient ascent [24]. Initialized near the true positions, the EM iterates converge to the global maximum; so this estimator, EM-GML, is a computable approximation of the GML estimator, the analogue of the UGIA-F benchmark drawing from the true positions. The UGIA-F estimator is realized by drawing position errors from the Gaussian distribution that achieves the inverse Fisher information matrix of the frame.

Three simulation studies follow. The first verifies Theorem 1 with a 2D Gaussian PSF, where summing $N$ frames grows the information and the GML estimator and UGIA-F both approach the CRB. The second verifies Theorem 2 with a 3D astigmatic PSF where the average intensity is swept and the same coincidence holds. The third examines the finite-photon regime, where the EM-GML estimator is necessarily biased, and shows that it attains an error below the unbiased CRB and below UGIA-F, indicating a gap between the GML estimate and the CRB that widens as the emitter density increases.

These results have two implications. First, EM-GML, used together with UGIA-F, provides a benchmark for practical localization algorithms, including neural networks that do not know the system parameters generating the frames. Second, the gap between the biased GML estimator and the unbiased CRB at high emitter density motivates the development of advanced localization algorithms, including neural networks that exploit this gap to achieve super spatiotemporal resolution.

**Notation.** Following Ref. [6], scalars are italic, e.g., $I$, $v(\boldsymbol{k})$, $\theta_{mj}$; vectors are bold italic, e.g., $\boldsymbol{\theta}$, $\boldsymbol{\theta}_m$, $\boldsymbol{u}$, $\boldsymbol{I}$, and the pixel index $\boldsymbol{k}$; matrices are bold upright capitals, e.g., $\mathbf{F}$, $\mathbf{I}_{dM}$; the pixel index set $\Omega$ is not bold, and the inverse Fisher information matrix is written $\mathbf{F}^{-1}$. A subscripted $\boldsymbol{k}$, e.g., $b_{\boldsymbol{k}}$, denotes a spatiotemporal density in photons/s/nm², and the argument form $(\boldsymbol{k})$, e.g., $b(\boldsymbol{k})$, denotes the corresponding photon count in pixel $\boldsymbol{k}$, that is, the density integrated over $\Delta_t$, $\Delta_x$, $\Delta_y$.

## 2. Imaging and noise model

**Camera and pixels.** A camera images $M$ point emitters on a pixel set $\Omega$ of $K = K_x \times K_y$ pixels with pixel index $\boldsymbol{k} = (k_x, k_y) \in \Omega \equiv \{0, 1, \dots, K_x - 1\} \times \{0, 1, \dots, K_y - 1\}$, pixel size $\Delta_x \times \Delta_y$ in nm$^2$, lateral field extents $L_x = K_x\Delta_x$ and $L_y = K_y\Delta_y$ in nm, and exposure time $\Delta_t$ in seconds.

**Positions, field of view, and parameter space.** Let $d \in \{2, 3\}$ be the spatial dimension. Emitter $m$ has position vector $\boldsymbol{\theta}_m = (\theta_{m1}, \dots, \theta_{md}) \in \mathbb{R}^d$, that is, $(x_m, y_m)$ in 2D and $(x_m, y_m, z_m)$ in 3D, and all positions are stacked into

$$\boldsymbol{\theta} = (\boldsymbol{\theta}_1^\top, \dots, \boldsymbol{\theta}_M^\top)^\top = (\theta_{11}, \dots, \theta_{1d}, \dots, \theta_{M1}, \dots, \theta_{Md})^\top \in \mathbb{R}^{dM}. \tag{1}$$

Let $\mathcal{F} \subset \mathbb{R}^d$ denote the *field of view (FOV)*, the region a single emitter can occupy, $\mathcal{F} = [0, L_x] \times [0, L_y]$ in 2D and $\mathcal{F} = [0, L_x] \times [0, L_y] \times [-L_z, L_z]$ in 3D where $2L_z$ is the axial range set by the PSF calibration. The *parameter space* is the set of $\delta$-separated configurations with every emitter in the FOV,

$$\Theta = \left\{\boldsymbol{\theta} \in \mathcal{F}^M\colon \min_{m \neq m'} \|\boldsymbol{\theta}_m - \boldsymbol{\theta}_{m'}\| \geq \delta\right\} \subset \mathbb{R}^{dM}. \tag{2}$$

The FOV $\mathcal{F}$ lives in $\mathbb{R}^d$, where one emitter sits; the parameter space $\Theta$ is the higher-dimensional set of full $M$-emitter configurations built from it.

**PSF and its pixel average.** Let $q_m(\boldsymbol{u})$ be the *point spread function (PSF)* of emitter $m$, the probability density that one of its detected photons lands at the camera-plane point $\boldsymbol{u} = (x, y)$, normalized to $\int q_m(\boldsymbol{u})\, \mathrm{d}\boldsymbol{u} = 1$, assuming a photon reaching the camera is detected with probability 1. The PSF is determined by the emitter position $\boldsymbol{\theta}_m$; in 3D the depth $z_m$ shapes it. The PSF may differ between emitters and be spatially variant. The detector integrates the PSF over each pixel; writing $p(\boldsymbol{k}) = [\Delta_x k_x, \Delta_x(k_x + 1)] \times [\Delta_y k_y, \Delta_y(k_y + 1)]$ for the area of pixel $\boldsymbol{k}$ and reusing the symbol $q_m$ by a mild abuse of notation, the *pixel-averaged PSF* is

$$q_m(\boldsymbol{k}) \;=\; \frac{1}{\Delta_x \Delta_y} \int_{p(\boldsymbol{k})} q_m(\boldsymbol{u})\, \mathrm{d}\boldsymbol{u}. \tag{3}$$

A working example for a 2D Gaussian PSF, its density and pixel-averaged form, is given in Appendix C where each $q_m(\boldsymbol{k})$ is $C^\infty$ in $\boldsymbol{\theta}$. $C^k$ denotes differentiability up to order $k$ throughout. The development below holds for a general PSF and is illustrated with the 2D Gaussian example.

**Signal photons.** The known emitter intensities are $I_1, \dots, I_M$ in photons/s, stacked as $\boldsymbol{I} = (I_1, \dots, I_M)^\top$. Following Ref. [6], the *average emitter intensity* is $I = M^{-1} \sum_{m=1}^M I_m$ and the fixed *intensity weights* are $\beta_m = I_m / I$ with $M^{-1} \sum_{m=1}^M \beta_m = 1$; set $Q(\boldsymbol{k}) = \sum_{m=1}^M \beta_m\, q_m(\boldsymbol{k})$ and $Q(\boldsymbol{u}) = \sum_{m=1}^M \beta_m\, q_m(\boldsymbol{u})$. The mean signal-photon density at camera-plane point $\boldsymbol{u}$ over the exposure is $s(\boldsymbol{u}) = \Delta_t\, I\, Q(\boldsymbol{u})$ in photons/nm². In pixel $\boldsymbol{k}$, the mean signal photon count is its integral over the pixel,

$$s(\boldsymbol{k}) \;=\; \int_{p(\boldsymbol{k})} s(\boldsymbol{u})\, \mathrm{d}\boldsymbol{u} \;=\; \Delta_t \Delta_x \Delta_y\, I\, Q(\boldsymbol{k}) \;=\; \Delta_t \Delta_x \Delta_y \sum_{m=1}^M I_m\, q_m(\boldsymbol{k}) \tag{4}$$

and the photon count in photons is $S(\boldsymbol{k}) \sim \mathrm{Poisson}(s(\boldsymbol{k}))$ independent across pixels where the pixel size $\Delta_x \Delta_y$ enters explicitly. The scale $I$ is the quantity sent to infinity in the high-intensity regime of Theorem 2 while the weights $\beta_m$, the PSF, and the noise parameters are held fixed.

**Background Poisson noise.** In pixel $\boldsymbol{k}$, let $b_{\boldsymbol{k}} > 0$, in photons/s/nm², be the spatiotemporal density of the background Poisson noise. The effect of spatially non-uniform background is denoted by the subscript $\boldsymbol{k}$. Its integral over the exposure and pixel area, $b(\boldsymbol{k}) = \Delta_t \Delta_x \Delta_y\, b_{\boldsymbol{k}}$ in photons, is the mean background count, and $B(\boldsymbol{k}) \sim \mathrm{Poisson}(b(\boldsymbol{k}))$ is the background Poisson noise sample.

**Gaussian readout noise.** Let $\mu_{\boldsymbol{k}}$ and $G_{\boldsymbol{k}}$, in photons/s/nm², be the spatiotemporal density mean and variance of the readout in pixel $\boldsymbol{k}$. The effect of pixel-variant readout noise in an sCMOS camera [25] is denoted by the subscript $\boldsymbol{k}$. The corresponding photon-count mean and variance are $\mu(\boldsymbol{k}) = \Delta_t \Delta_x \Delta_y\, \mu_{\boldsymbol{k}}$ and $G(\boldsymbol{k}) = \Delta_t \Delta_x \Delta_y\, G_{\boldsymbol{k}}$ in photons, and the readout Gaussian noise sample is $W(\boldsymbol{k}) \sim \mathcal{N}(\mu(\boldsymbol{k}), G(\boldsymbol{k}))$, independent of the signal and background. The spatially invariant model in Refs. [5, 6] is the special case for $b_{\boldsymbol{k}} = b$, $\mu_{\boldsymbol{k}} = \mu$, $G_{\boldsymbol{k}} = G$.

**Mixed Poisson–Gaussian noise and its Poisson approximation.** Adjust the readout mean to equal its variance, $\mu(\boldsymbol{k}) = G(\boldsymbol{k})$, equivalently $\mu_{\boldsymbol{k}} = G_{\boldsymbol{k}}$; in practice the known readout offset is subtracted from the recorded frame and its variance added back. The combined background-plus-readout noise then has *equal mean and variance*, and when the per-pixel photon count is not small (e.g., $b(\boldsymbol{k}) \geq 10$ photons), is well approximated by a Poisson variable [5, 6]. With a mild abuse of notation we retain the symbols $b_{\boldsymbol{k}}$, $b(\boldsymbol{k})$, $B(\boldsymbol{k})$ for this *mixed* noise, whose spatiotemporal density mean becomes $b_{\boldsymbol{k}} + G_{\boldsymbol{k}}$, written again as $b_{\boldsymbol{k}}$, so the mixed-noise mean count is $b(\boldsymbol{k}) = \Delta_t \Delta_x \Delta_y\, b_{\boldsymbol{k}}$, that is, the original background count plus $G(\boldsymbol{k})$, and the mixed-noise sample is $B(\boldsymbol{k}) \sim \text{Poisson}\big(b(\boldsymbol{k})\big)$ in photons.

**Total frame (working model).** The pixel value is the signal plus the mixed noise,

$$V(\boldsymbol{k}) = S(\boldsymbol{k}) + B(\boldsymbol{k}) \sim \text{Poisson}\big(v(\boldsymbol{k})\big) \tag{5}$$

with $v(\boldsymbol{k}) = s(\boldsymbol{k}) + b(\boldsymbol{k})$, independent across pixels. No correction symbol is needed since $v(\boldsymbol{k})$ and $V(\boldsymbol{k})$ already denote the signal-plus-noise mean and sample. This Poisson model is the one the estimator maximizes and under which the Fisher information $\mathbf{F}$ is derived; Section 3.5 quantifies the vanishing error relative to the true Poisson–Gaussian frame.

**Per-frame log-likelihood, score, and Fisher information.** By the pixel independence of Eq. (5), the frame $V$ has probability mass function

$$f_V(V;\boldsymbol{\theta}) = \prod_{\boldsymbol{k}\in\Omega} \frac{v(\boldsymbol{k})^{V(\boldsymbol{k})}\, e^{-v(\boldsymbol{k})}}{V(\boldsymbol{k})!}, \tag{6}$$

which is also the likelihood function of $\boldsymbol{\theta}$. Dropping the $\boldsymbol{\theta}$-independent constant as in all log-likelihood functions henceforth, the per-frame log-likelihood, score, and Fisher information are

$$\ell_1(\boldsymbol{\theta}) \;=\; \sum_{\boldsymbol{k}\in\Omega} [V(\boldsymbol{k}) \ln v\,(\boldsymbol{k}) - v(\boldsymbol{k})]\,, \tag{7}$$

$$\boldsymbol{s}_1(\boldsymbol{\theta}) \;=\; \nabla_{\boldsymbol{\theta}}\, \ell_1(\boldsymbol{\theta}) \;=\; \sum_{\boldsymbol{k}\in\Omega} \left(\frac{V(\boldsymbol{k})}{v(\boldsymbol{k})} - 1\right) \frac{\partial v(\boldsymbol{k})}{\partial \boldsymbol{\theta}}, \tag{8}$$

$$\mathbf{F}(\boldsymbol{\theta}) \;=\; \mathbb{E}_{\boldsymbol{\theta}}[\boldsymbol{s}_1 \boldsymbol{s}_1^{\top}] \;=\; \sum_{\boldsymbol{k}\in\Omega} \frac{1}{v(\boldsymbol{k})} \frac{\partial v(\boldsymbol{k})}{\partial \boldsymbol{\theta}} \frac{\partial v(\boldsymbol{k})}{\partial \boldsymbol{\theta}}^{\top} \;\in\; \mathbb{R}^{dM\times dM}, \tag{9}$$

using $\mathbb{E}[V(\boldsymbol{k})] = \text{Var}[V(\boldsymbol{k})] = v(\boldsymbol{k})$, pixel independence, and $\partial v(\boldsymbol{k})/\partial\boldsymbol{\theta} = \partial s(\boldsymbol{k})/\partial\boldsymbol{\theta} = \Delta_t\Delta_x\Delta_y \sum_m I_m\; \partial q_m(\boldsymbol{k})/\partial\boldsymbol{\theta}$ since the noise mean $b(\boldsymbol{k})$ is constant in $\boldsymbol{\theta}$. Eq. (9) is the closed-form Fisher information matrix (FIM) $\mathbf{F}$ derived for SMLM in Ref. [5]. Because the intensities are known, it is case (i) of Ref. [6] where the full FIM reduces to its position block, $\mathbf{F} = \mathbf{F}_{11}$, and the per-pixel weight $1/v(\boldsymbol{k})$ coincides with its Corollary 2 where the mixed-noise mean $b(\boldsymbol{k})$ is in the role of background-plus-readout.

**Cramér–Rao bound and lower bound on RMSE.** We use the following objects. The *Cramér–Rao bound (CRB)* [13, 14] is the *matrix* $\mathbf{F}(\boldsymbol{\theta})^{-1}$, so the covariance matrix of any unbiased estimator $\hat{\boldsymbol{\theta}}$ is bounded below in the Loewner ordering, $\text{Cov}\big(\hat{\boldsymbol{\theta}}\big) \succeq \mathbf{F}^{-1}$. The *Cramér–Rao lower bound (CRLB)* of a single coordinate $\theta_{mj}$ is the standard deviation

$$\Delta\big(\theta_{mj}\big) \;=\; \big([\mathbf{F}^{-1}]_{mj,mj}\big)^{1/2}, \tag{10}$$

the square root of the corresponding diagonal element of $\mathbf{F}^{-1}$, which already accounts for coupling to the other coordinates and is in general far larger than $[\mathbf{F}]_{mj,mj}^{-1/2}$.

Because the emitter indices $m$ of the estimated position $\hat{\boldsymbol{\theta}}_m$ and the true position $\boldsymbol{\theta}_m$ are known one-to-one for all the estimators to be investigated, the accuracy of an estimator $\hat{\boldsymbol{\theta}}$ is scored by the *root-mean-square error (RMSE)* $R$ in nm, defined by $R^2 = \mathbb{E}\left(M^{-1}\sum_{m=1}^{M}\|\hat{\boldsymbol{\theta}}_m - \boldsymbol{\theta}_{0,m}\|^2\right)$. The mean $\mathbb{E}$ over the photon count of signal and noise is obtained by a Monte-Carlo average over frames in practice.

For any *unbiased* estimator each coordinate variance is bounded below by its CRLB, $\mathrm{Var}(\hat{\theta}_{mj}) \geq \Delta(\theta_{mj})^2$, giving the *lower bound on the RMSE*

$$R^2 = \frac{1}{M}\sum_{m=1}^{M}\sum_{j=1}^{d}\mathrm{Var}(\hat{\theta}_{mj}) \;\geq\; \frac{1}{M}\sum_{m=1}^{M}\sum_{j=1}^{d}\Delta\left(\theta_{mj}\right)^2 = \frac{1}{M}\mathrm{tr}\;(\mathbf{F}^{-1}) \tag{11}$$

with equality for an efficient estimator. UGIA-F attains it by construction, so $R^2_{\text{UGIA-F}} = M^{-1}\,\mathrm{tr}(\mathbf{F}^{-1})$ is the *information-achieving lower benchmark* for unbiased estimators. The matrix CRB $\mathbf{F}(\boldsymbol{\theta})^{-1}$ is the fundamental object since it carries the inter-emitter correlations that the diagonal CRLB and the trace omit. The efficiency statements below refer to it.

**Total Fisher information.** For $N$ independent frames at the same $\boldsymbol{\theta}$, the total log-likelihood is $\ell_N(\boldsymbol{\theta}) = \sum_{n=1}^{N}\ell_1^{(n)}(\boldsymbol{\theta})$, the total score is $\boldsymbol{s}_N(\boldsymbol{\theta}) = \nabla_{\boldsymbol{\theta}}\ell_N(\boldsymbol{\theta}) = \sum_{n=1}^{N}\boldsymbol{s}_1^{(n)}(\boldsymbol{\theta})$, and the Fisher informations add. That is, the *total Fisher information* is $\mathbf{F}_{\text{total}} = N\,\mathbf{F}(\boldsymbol{\theta})$, more generally $N\,\mathbf{F}(I)$ at average intensity $I$ as in Section 3.3. Throughout, "*the information grows*" means the minimum eigenvalue $\lambda_{\min}(\mathbf{F}_{\text{total}}) \to \infty$, equivalently the CRB vanishes, $\mathbf{F}_{\text{total}}^{-1} \to \mathbf{0}$, driven by $N \to \infty$ and/or $I \to \infty$.

## 3. Theory

### 3.1 Estimators

The following estimators are considered.

**GML** (global maximum likelihood): $\hat{\boldsymbol{\theta}}_N^{\text{GML}} \in \arg\max_{\boldsymbol{\theta}\in\Theta}\ell_N(\boldsymbol{\theta})$, the global maximizer.

**LML** (local maximum likelihood): any stationary point reached by an iterative maximizer from a generic initialization, which is unnecessary to be global.

**UGIA-F** (unbiased Gaussian information-achieving estimator for a frame): by construction

$$\hat{\boldsymbol{\theta}}^F \;=\; \boldsymbol{\theta}_0 + \mathbf{F}(\boldsymbol{\theta}_0)^{-1/2}\,\boldsymbol{z} \tag{12}$$

where $\boldsymbol{z} \sim \mathcal{N}(\mathbf{0}, \mathbf{I}_{dM})$, hence $\hat{\boldsymbol{\theta}}^F \sim \mathcal{N}(\boldsymbol{\theta}_0, \mathbf{F}(\boldsymbol{\theta}_0)^{-1})$, an *unbiased* estimator whose covariance equals the per-frame CRB. It was introduced in Ref. [7], where it is constructed from the singular value decomposition (SVD) of the FIM, and is a theoretical benchmark and not deployable as it uses $\boldsymbol{\theta}_0$.

### 3.2 Regularity conditions

Fix the true positions $\boldsymbol{\theta}_0$. The following conditions are assumed at $\boldsymbol{\theta}_0$.

**C1. Parameter space.** $\Theta$, defined in Section 2, is compact and $\boldsymbol{\theta}_0 \in \mathrm{int}\,\Theta$, i.e., every emitter is strictly inside the FOV $\mathcal{F}$ and all pairwise separations strictly exceed $\delta$.

**C2. Positivity and smoothness.** There exists $v_{\min} > 0$ with $v(\boldsymbol{k}) \geq v_{\min}$ for all $\boldsymbol{k}$ and all $\boldsymbol{\theta} \in \Theta$; and $v(\boldsymbol{k}) \in C^3(\Theta)$, or more weakly the model is differentiable in quadratic mean (DQM), which suffices for the asymptotics and admits $C^2$ spline PSFs, see Section 3.4.

**C3. Information.** $\mathbf{F}(\boldsymbol{\theta}_0) \succ \mathbf{0}$, positive definite.

**C4. Uniform regularity.** The third-order partial derivatives of $\ell_1(\boldsymbol{\theta})$ are bounded by an integrable envelope uniformly on a neighborhood of $\boldsymbol{\theta}_0$.

Section 3.4 verifies that C1–C4 hold for the SMLM Poisson model.

### 3.3 Main results

The main results of this paper are as follows.

*Theorem 1* (Asymptotic efficiency, $N \to \infty$, multiple frames): Under (C1)–(C4), given $N$ i.i.d. frames from the Poisson working model at $\boldsymbol{\theta}_0$, the GML estimator $\hat{\boldsymbol{\theta}}_N^{\mathrm{GML}}$ satisfies as $N \to \infty$,
(i) $\hat{\boldsymbol{\theta}}_N^{\mathrm{GML}} \xrightarrow{p} \boldsymbol{\theta}_0$;
(ii) $\sqrt{N}\left(\hat{\boldsymbol{\theta}}_N^{\mathrm{GML}} - \boldsymbol{\theta}_0\right) \xrightarrow{d} \mathcal{N}(\mathbf{0},\ \mathbf{F}(\boldsymbol{\theta}_0)^{-1})$;
(iii) $\mathbb{E}\left(\hat{\boldsymbol{\theta}}_N^{\mathrm{GML}}\right) \to \boldsymbol{\theta}_0$ and bias $= O(N^{-1})$.
Equivalently, for large $N$

$$\hat{\boldsymbol{\theta}}_N^{\mathrm{GML}} \approx \mathcal{N}\left(\boldsymbol{\theta}_0,\ \left(N\mathbf{F}(\boldsymbol{\theta}_0)\right)^{-1}\right), \tag{13}$$

which is asymptotically Gaussian and unbiased, attaining the CRB $\mathbf{F}_{\text{total}}^{-1} = (N\mathbf{F})^{-1}$. □

The three GML properties (i)-(iii) are consistency, asymptotic normality and efficiency, and asymptotic unbiasedness, respectively.

*Theorem 2* (Asymptotic efficiency, $I \to \infty$, single frame): Let $I \to \infty$ with the weights $\beta_m$ fixed. Then $\mathbf{F}(I) \succ \mathbf{0}$ with $\lambda_{\min}\left(\mathbf{F}(I)\right) \to \infty$, equivalently $\mathbf{F}(I)^{-1} \to \mathbf{0}$, and the single-frame GML estimator satisfies

$$\mathbf{F}(I)^{1/2}\left(\hat{\boldsymbol{\theta}}^{\mathrm{GML}} - \boldsymbol{\theta}_0\right) \xrightarrow{d} \mathcal{N}(\mathbf{0}, \mathbf{I}_{dM}), \tag{14}$$

i.e., $\hat{\boldsymbol{\theta}}^{\mathrm{GML}} \approx \mathcal{N}(\boldsymbol{\theta}_0, \mathbf{F}(I)^{-1})$, attaining the CRB. □

Combining the two limits of Theorems 1 and 2, the total information is $N\mathbf{F}(I)$ with $N$ frames at average intensity $I$. The proofs of Theorems 1 and 2 are given in Appendix A.

UGIA-F is a draw from $\mathcal{N}(\boldsymbol{\theta}_0, \mathbf{F}(\boldsymbol{\theta}_0)^{-1})$ by construction of Eq. (12). By Theorem 1 (ii) and Theorem 2, the GML estimator converges in distribution to the same law. Therefore, the following corollary is obtained.

*Corollary 1*: As $N \to \infty$ and/or $I \to \infty$,

$$\hat{\boldsymbol{\theta}}^{\mathrm{GML}} \ \xrightarrow{d}\ \mathcal{N}(\boldsymbol{\theta}_0, \mathbf{F}^{-1}) \ \stackrel{d}{=}\ \hat{\boldsymbol{\theta}}^{F}, \tag{15}$$

the GML estimator and the UGIA-F estimator converge to the identical distribution. □

### *3.4 Verification of conditions C1–C4*

The four conditions are verified for the SMLM Poisson model.

**C1. Parameter space.** Imposed by construction, restricting to the interior and $\delta$-separated configurations. Compactness and the interior truth follow.

**C2. Positivity and smoothness.** Because the mixed-noise mean is strictly positive, the total mean is bounded below,

$$v(\boldsymbol{k}) = s(\boldsymbol{k}) + b(\boldsymbol{k}) \ \geq\ b(\boldsymbol{k}) \ \geq\ v_{\min} > 0 \tag{16}$$

for all $\boldsymbol{\theta}, \boldsymbol{k}$ with $v_{\min} = \min_{\boldsymbol{k}\in\Omega} b(\boldsymbol{k}) > 0$, a minimum over the finite pixel set $\Omega$ of strictly positive mixed-noise means, each $b(\boldsymbol{k}) = \Delta_t\Delta_x\Delta_y\, b_{\boldsymbol{k}} > 0$ . So, the Poisson mean never approaches zero, thereby eliminating the degeneracy that would otherwise spoil regularity at low counts. Regarding smoothness, for an *analytic PSF* such as Gaussian, scalar-diffraction/Airy [26], Gibson–Lanni [27], vectorial [28], astigmatic [19], double-helix [20], or tetrapod [21], each $q_m(\boldsymbol{k})$ and hence $v(\boldsymbol{k})$ is $C^\infty$ in $\boldsymbol{\theta}$, so C2 and C4 hold with room to spare. A *cubic-spline experimental PSF* [22] gives only $v(\boldsymbol{k}) \in C^2$; this is still covered because the asymptotics require only differentiability in quadratic mean by following Refs. [15, 17], which is weaker than the classical $C^3$ condition of Ref. [14]. Therefore, smoothness is available for these PSFs.

**C3. Information.** From the closed-form FIM in Eq. (9),

$$\mathbf{F}(\boldsymbol{\theta}_0) = \sum_{\boldsymbol{k}\in\Omega} \frac{1}{v(\boldsymbol{k})} \frac{\partial v(\boldsymbol{k})}{\partial \boldsymbol{\theta}} \frac{\partial v(\boldsymbol{k})}{\partial \boldsymbol{\theta}}^{\top} \succ \mathbf{0} \Leftrightarrow \left\{ \frac{\partial v(\boldsymbol{k})}{\partial \boldsymbol{\theta}} (\theta_0) \right\}_{\boldsymbol{k}\in\Omega} \text{ span } \mathbb{R}^{dM}. \quad (17)$$

For $\delta$-separated emitters the $dM$ position-gradient fields have distinct, essentially disjoint supports and are linearly independent as $K$-vectors, so $\mathbf{F} \succ \mathbf{0}$. This is directly *checkable numerically* from the coded FIM with the smallest eigenvalue greater than zero; $\lambda_{\min}(\mathbf{F}) \to 0$ and the CRB $\mathbf{F}^{-1}$ blows up as $\delta \to 0$, a loss of resolution. *In 3D*, where $d = 3$, the condition additionally requires the PSF to be *axially informative*, that is, $\partial v(\boldsymbol{k})/\partial z_m$ is nonzero and linearly independent of the lateral gradients, which holds only for $z$-encoding PSFs such as astigmatic, double-helix, or tetrapod. Note that an axially symmetric PSF at focus carries no axial information, so $\mathbf{F}$ is singular in $z$ and efficiency fails for the axial coordinate while lateral efficiency is unaffected.

**C4. Uniform regularity.** It follows from Eq. (7) that the third derivative $\nabla^3 \ell_1$ is a finite sum of terms built from $V(\boldsymbol{k})/v(\boldsymbol{k})^j$, $j \le 3$, and the derivatives of $v(\boldsymbol{k})$ up to the third order. With $v(\boldsymbol{k}) \ge v_{\min}$ (C2), bounded PSF derivatives on compact $\Theta$, and $\mathbb{E}[V(\boldsymbol{k})] = v(\boldsymbol{k}) < \infty$, these are dominated by an integrable envelope, giving the required control of the Taylor remainder.

### *3.5 Remark: the Gaussian-to-Poisson approximation*

Theorems 1 and 2 are exact under the Poisson working model. When the data are generated from the *true* Poisson–Gaussian frame of signal plus Poisson background plus Gaussian readout but estimated with the Poisson likelihood, the GML estimator is a *quasi-MLE* [29], which converges to the pseudo-true value $\boldsymbol{\theta}_*$ minimizing the KL divergence, i.e., the relative entropy, from the true law where the sandwich covariance is $\mathbf{A}^{-1}\mathbf{B}\,\mathbf{A}^{-1}$ with $\mathbf{A} = -\mathbb{E}(\nabla^2 \ell_1)$ and $\mathbf{B} = \mathrm{Cov}(\boldsymbol{s}_1)$. The mean-adjustment $\mu(\boldsymbol{k}) = G(\boldsymbol{k})$ makes the first two moments of $V(\boldsymbol{k})$ match the Poisson working model at $\boldsymbol{\theta}_0$, i.e., both mean and variance $v(\boldsymbol{k})$; hence $\boldsymbol{\theta}_* = \boldsymbol{\theta}_0$ and $\mathbf{A} = \mathbf{B} = \mathbf{F}$ to this order, so the sandwich collapses to the CRB $\mathbf{F}^{-1}$ and the theorem holds. The residual mismatch enters only through third and higher moments and is $O(1/I)$, which vanishes as $I \to \infty$ in the signal-dominated regime. Consequently, (i) on data simulated from the Poisson working model the theorem holds *exactly*, recommended for the clean verification; (ii) on the true Poisson–Gaussian frame it holds *up to* $O(1/I)$, which the simulation can quantify.

### *3.6 The gap between GML and UGIA-F with finite information*

While Theorems 1 and 2 and Corollary 1 present the asymptotic properties of GML and UGIA-F and their relationship in the infinite-photon regime, it is interesting in practice to know the gap between GML and UGIA in the finite-photon regime. We hypothesize below.

*Hypothesis 1*: In the finite-information regime where $N$ and $I$ are finite, the GML estimator is *biased* and attains a smaller RMSE than the unbiased UGIA-F estimator. Their deviation grows with the overlap severity of PSFs. □

As indicated by Hypothesis 1, the finite-information regime is where GML and UGIA-F differ. There $\hat{\boldsymbol{\theta}}^{\mathrm{GML}}$ is *biased*, and because the CRB constrains only *unbiased* estimators, the biased GML can attain a smaller error than the unbiased $\hat{\boldsymbol{\theta}}^{F}$, increasingly so as the CRB inflates with growing overlap of PSFs, $\mathrm{tr}(\mathbf{F}^{-1}) \to \infty$, equivalently $\lambda_{\min}(\mathbf{F}) \to 0$. Since $R^2_{\text{UGIA-F}} = \mathrm{tr}(\mathbf{F}^{-1})/M$ in Eq. (11), this is the blow-up of the unbiased UGIA-F. This deviation from the common asymptotic limit, which grows with density, is exactly *Hypothesis 1* that GML outperforms UGIA-F at high density. It is proved to shrink to zero as the total Fisher information grows, $\lambda_{\min}\big(N\,\mathbf{F}(I)\big) \to \infty$, by Theorems 1 and 2.

## 4. The EM algorithm

We shall numerically verify Theorems 1 and 2, Corollary 1, and Hypothesis 1. The GML estimator of Section 3.1 maximizes the frame log-likelihood $\ell_1(\boldsymbol{\theta})$ of Eq. (7) jointly over the $dM$ position coordinates, which is computationally difficult at high emitter density. The expectation-maximization (EM) algorithm [23, 24] exploits the Poisson properties and the mutual independence of the per-emitter photon emissions in the frame model of Section 2, which replaces this joint maximization by a sequence of single-emitter localizations with monotone likelihood ascent.

Write the mean per-emitter signal count $s_m(\boldsymbol{k}) = \Delta_t \Delta_x \Delta_y\, I_m q_m(\boldsymbol{k})$ so that $S_m(\boldsymbol{k}) \sim \text{Poisson}\big(s_m(\boldsymbol{k})\big)$ is the image of emitter $m$ alone, $s(\boldsymbol{k}) = \sum_{m=1}^{M} s_m(\boldsymbol{k})$, and the frame is $V(\boldsymbol{k}) = \sum_{m=1}^{M} S_m(\boldsymbol{k}) + B(\boldsymbol{k})$ with $B(\boldsymbol{k}) \sim \text{Poisson}\big(b(\boldsymbol{k})\big)$ the mixed noise of Section 1. The variables $S_m(\boldsymbol{k})$ and $B(\boldsymbol{k})$ are mutually independent across emitters and pixels.

### 4.1 Decomposition and ascent

Since the per-emitter images $S_m(\boldsymbol{k})$ are mutually independent, their joint probability distribution can be written as

$$\prod_{m=1}^{M} f_{S_m}(S_m; \boldsymbol{\theta}_m) = \prod_{m=1}^{M} f_{S_m|V}(S_m; \boldsymbol{\theta}_m)\, f_V(V; \boldsymbol{\theta}) \tag{18}$$

where $f_{S_m}(S_m; \boldsymbol{\theta}_m) = \prod_{\boldsymbol{k}\in\Omega} s_m(\boldsymbol{k})^{S_m(\boldsymbol{k})} e^{-s_m(\boldsymbol{k})} / S_m(\boldsymbol{k})!$ is the probability mass function of the single-emitter noise-free image $S_m(\boldsymbol{k})$ and $f_{S_m|V}(S_m; \boldsymbol{\theta}_m)$ is its probability mass function conditioned on the frame $V(\boldsymbol{k})$, $\boldsymbol{k} \in \Omega$. Taking the logarithm of both sides and rearranging, with $\ell_1(\boldsymbol{\theta}) = \ln f_V(V; \boldsymbol{\theta})$, yields the frame log-likelihood decomposition

$$\ell_1(\boldsymbol{\theta}) = \sum_{m=1}^{M} \ln f_{S_m}(S_m; \boldsymbol{\theta}_m) - \sum_{m=1}^{M} \ln f_{S_m|V}(S_m; \boldsymbol{\theta}_m). \tag{19}$$

Taking the conditional expectation of both sides given the frame and the current estimate $\boldsymbol{\theta}^{(l)}$, and using that $\ell_1(\boldsymbol{\theta})$ is fixed given the frame, yields

$$\ell_1(\boldsymbol{\theta}) = Q\big(\boldsymbol{\theta}; \boldsymbol{\theta}^{(l)}\big) - \sum_{m=1}^{M} \mathbb{E}\big[\ln f_{S_m|V}(S_m; \boldsymbol{\theta}_m)\,\big|V, \boldsymbol{\theta}^{(l)}\big] \tag{20}$$

where

$$Q\big(\boldsymbol{\theta}; \boldsymbol{\theta}^{(l)}\big) = \sum_{m=1}^{M} \mathbb{E}\big[\ln f_{S_m}(S_m; \boldsymbol{\theta}_m)\big|V, \boldsymbol{\theta}^{(l)}\big] \tag{21}$$

is the EM auxiliary function. By Jensen's inequality the second sum is maximized over $\boldsymbol{\theta}$ at $\boldsymbol{\theta} = \boldsymbol{\theta}^{(l)}$, so any $\boldsymbol{\theta}^{(l+1)}$ that increases $Q$ does not decrease $\ell_1$, that is, $\ell_1\big(\boldsymbol{\theta}^{(l+1)}\big) \geq \ell_1\big(\boldsymbol{\theta}^{(l)}\big)$. The auxiliary function separates across emitters, so its maximization reduces to $M$ independent single-emitter problems for $m = 1, \dots, M$,

$$\hat{\boldsymbol{\theta}}_m^{(l+1)} = \arg\max_{\boldsymbol{\theta}_m} \mathbb{E}\big[\ln f_{S_m}(S_m; \boldsymbol{\theta}_m)\big|V, \boldsymbol{\theta}^{(l)}\big]. \tag{22}$$

### 4.2 E-step and M-step

The E-step evaluates the conditional expectation. Since $S_m(\boldsymbol{k})$ and $B(\boldsymbol{k})$ are mutually independent Poisson variables whose sum is $V(\boldsymbol{k})$, by a Poisson property [30] the conditional mean of $S_m(\boldsymbol{k})$ given the frame is,

$$\hat{S}_m^{(l+1)}(\boldsymbol{k}) = \mathbb{E}\big[S_m(\boldsymbol{k})|V, \boldsymbol{\theta}^{(l)}\big] = \frac{s_m^{(l)}(\boldsymbol{k})}{v^{(l)}(\boldsymbol{k})} V(\boldsymbol{k}) = \frac{I_m q_m^{(l)}(\boldsymbol{k})}{\sum_{i=1}^{M} I_i\, q_i^{(l)}(\boldsymbol{k}) + b_{\boldsymbol{k}}} V(\boldsymbol{k}) \tag{23}$$

where the superscript $(l)$ denotes evaluation at $\boldsymbol{\theta}^{(l)}$, $v^{(l)}(\boldsymbol{k}) = \sum_{i=1}^{M} s_i^{(l)}(\boldsymbol{k}) + b(\boldsymbol{k})$, and $b_{\boldsymbol{k}}$ is the spatiotemporal density mean of mixed noise in Section 1. The fraction is the share of the frame count attributed to emitter $m$, so the E-step assigns to each emitter a noise-free estimated image.

Since $\ln f_{S_m}(S_m; \boldsymbol{\theta}_m)$ is linear in $S_m(\boldsymbol{k})$, its conditional expectation replaces $S_m(\boldsymbol{k})$ by $\hat{S}_m^{(l+1)}(\boldsymbol{k})$. The M-step then maximizes the single-emitter likelihood of the estimated image, which gives the noise-free single-emitter ML estimate

$$\hat{\boldsymbol{\theta}}_m^{(l+1)} = \arg\max_{\boldsymbol{\theta}_m} \sum_{\boldsymbol{k}\in\Omega} \big[\hat{S}_m^{(l+1)}(\boldsymbol{k}) \ln s_m(\boldsymbol{k}) - s_m(\boldsymbol{k})\big]. \tag{24}$$

The two single-emitter ML solvers by the gradient ascent and Newton are given in Appendix B.

### 4.3 Algorithm

Given an initial estimate $\boldsymbol{\theta}^{(0)}$, repeat two steps. The E-step computes $\hat{S}_m^{(l+1)}(\boldsymbol{k})$ for all $m$ and $\boldsymbol{k} \in \Omega$ by Eq. (23) from the frame $V(\boldsymbol{k})$ and the current $\boldsymbol{\theta}^{(l)}$. The M-step computes $\hat{\boldsymbol{\theta}}_m^{(l+1)}$ for all $m$ by Eq. (24) from the $\hat{S}_m^{(l+1)}(\boldsymbol{k})$. Terminate when $\|\boldsymbol{\theta}^{(l+1)} - \boldsymbol{\theta}^{(l)}\| < \varepsilon$ or $l \geq l^*$ for prespecified $\varepsilon > 0$ and $l^* > 0$.

### 4.4 Properties

The EM algorithm has three properties of interest here. First, as a likelihood-ascent search it converges in general to an LML estimate called EM-LML; but if the initial estimate $\boldsymbol{\theta}^{(0)}$ is close to the true positions $\boldsymbol{\theta}_0$, it converges to the GML estimate, called EM-GML, in the high-information regime, which is the regime verified in Section 5. Second, the likelihood is monotonically non-decreasing along the iterates by Eq. (22), which improves the localization from step to step and guarantees convergence. Third, the EM algorithm is computationally simple and fast, since Eqs. (23) and (24) decompose the joint $M$-emitter localization into $M$ single-emitter localizations consisting of the E-step a pixelwise fraction and the M-step a noise-free single-emitter ML. The compute for all emitters can be run in parallel as a batch and be efficiently implementable on a GPU.

By targeting the global maximum, EM-GML is the operational realization of the asymptotically efficient estimator of Theorems 1 and 2 and provides the general accuracy bound against biased, unbiased, and learned estimators. From an arbitrary initialization the EM converges in general to a local maximum, EM-LML, which is biased and does not carry the efficiency guarantee, so the near-true initialization is used throughout Section 5.

## 5. Simulation verification

### 5.1 Properties to be verified

Fix a high-density configuration, e.g., $M = 80$ near the resolution limit. Sweep the information, i.e., raise $\lambda_{\min}(\mathbf{F}_{\text{total}})$ by increasing the number of frames $N$ or the average intensity $I$. For each setting, draw $L$ Monte-Carlo replicates and estimate positions with EM initialized near the truth, so the *global* maximum is found by EM-GML. The following five properties shall be verified by simulations.

1. **Consistency and unbiasedness.** The RMSE $\to 0$ and the bias $\left\|\bar{\hat{\boldsymbol{\theta}}} - \boldsymbol{\theta}_0\right\|$ stays at or below the Monte-Carlo floor RMSE$/\sqrt{L}$, so it is negligible relative to the standard deviation. The

$O(1/N)$ bias rate, below the resolution of $L$ replicates, is a property of the GML by Theorem 1.

2. **Efficiency.** The ratio of the empirical per-coordinate variance to the squared CRLB $\Delta(\theta_{mj})^2$, Eq. (10), tends to 1, equivalently the whitened residual of property 3 has unit per-coordinate variance. The full covariance equality $\mathrm{Cov}(\hat{\boldsymbol{\theta}}) \to (N\mathbf{F})^{-1}$ follows from Theorem 1 and is not separately estimated, since the full $dM \times dM$ covariance is beyond the Monte-Carlo replicate count, while $\mathbf{F}$ is near block-diagonal for resolvable emitters so the per-coordinate variance is the informative part.
3. **Normality.** A QQ-plot tests the normality of the whitened residuals $\mathbf{z} = \mathbf{F}^{1/2}(\hat{\boldsymbol{\theta}} - \boldsymbol{\theta}_0)$, pooled over the emitters and the $L$ replicates, which are $\mathcal{N}(\mathbf{0}, \mathbf{I})$ under Theorems 1 and 2. The same residual verifies both the normality here and the per-coordinate efficiency of property 2.
4. **GML ≡ UGIA-F.** $R(\hat{\boldsymbol{\theta}}^{\mathrm{GML}})$ and $R(\hat{\boldsymbol{\theta}}^{F})$ converge as the information grows, the coincidence of the biased GML and the unbiased benchmark in the high-information limit by Corollary 1.
5. **GML outperforms UGIA-F at high density.** In the finite-information regime of finite $N$ and $I$ the biased EM-GML attains a smaller RMSE than the unbiased UGIA-F, so $R(\hat{\boldsymbol{\theta}}^{\mathrm{GML}}) < R(\hat{\boldsymbol{\theta}}^{F})$ by Hypothesis 1, and the gap widens as the emitter density increases.

EM-GML with near-true initialization is the appropriate global-maximum finder for this verification, since efficiency is a property of the global maximizer only; EM from random initialization converges to a local maximum, i.e., EM-LML, and is not used.

The gap of Hypothesis 1 is measured by the RMSE under the known one-to-one correspondence between the estimated and true positions, so each estimator is compared to the CRB on the same footing. Each configuration is summarized by its signal-to-noise ratio (SNR), computed from the system parameters following the information-sufficient SNR of Ref. [31] for the 2D PSF and its 3D astigmatic extension by Ref. [32].

### *5.2 Approach*

**EM-GML.** The EM-GML estimator is the EM algorithm initialized with a small jitter around the true positions, each emitter initially placed within $\delta/2$ nm of its own true position, so the global maximum is found. The EM algorithm can also be initialized by an emitter index-agnostic estimator such as successive interference cancellation (SIC) of Ref. [24] to reach the global maximum in practice; so this near-true initialization stands in for an achievable one that is the EM-GML analogue of UGIA-F drawing from the true positions.

**UGIA-F.** The UGIA-F estimator generates the sample according to $\hat{\boldsymbol{\theta}}^{F} \sim \mathcal{N}(\boldsymbol{\theta}_0, \mathbf{F}(\boldsymbol{\theta}_0)^{-1})$ where $\mathbf{F}(\boldsymbol{\theta}_0)$ is the FIM of a frame.

**RMSE.** The estimated and true positions are in known one-to-one correspondence, so the accuracy of an estimator is scored by the RMSE

$$\mathrm{RMSE}(\hat{\boldsymbol{\theta}}, \boldsymbol{\theta}_0)^2 = \frac{1}{M}\sum_{m=1}^{M} \left\| \hat{\boldsymbol{\theta}}_m - \boldsymbol{\theta}_{0,m} \right\|^2, \tag{25}$$

and the reported value is the root mean of Eq. (25) over the $L$ Monte-Carlo replicates.

**Information-sufficient SNR.** The information-sufficient SNR is proposed and derived for the 2D Gaussian PSF [31] and the 3D astigmatic PSF [32] for the uniform noise and the equal emitter intensity. We extend the SNR to the non-uniform noise and emitter-dependent intensity of this model. The signal quality of each configuration is quantified by a per-emitter information-sufficient SNR.

The Poisson background and Gaussian readout enter through their summed variance density $b_{\boldsymbol{k}} + G_{\boldsymbol{k}}$, and the local noise density seen by emitter $m$ is the PSF-weighted mean

$$\bar{n}_m = \frac{\sum_{\boldsymbol{k}\in\Omega} q_m(\boldsymbol{k})\,(b_{\boldsymbol{k}} + G_{\boldsymbol{k}})}{\sum_{\boldsymbol{k}\in\Omega} q_m(\boldsymbol{k})}. \tag{26}$$

The per-emitter intensity-to-noise ratio is $\gamma_m = I_m/\bar{n}_m$, the figure of merit of the PSF at the emitter is $\eta_m$, and the *per-emitter* information-sufficient SNR is defined as $v_m = \eta_m\gamma_m$ with $\mathrm{SNR}_m = 10\log_{10} v_m$ in dB.

For the 2D Gaussian PSF of Appendix C with width $\sigma$ the figure of merit [31] is

$$\eta_m = \frac{\eta_s}{\sigma^2} \tag{27}$$

where $\eta_s = -\rho/[2\pi\ln(1-\rho)] = 0.0791$ is the figure of merit of the unit-variance Gaussian at the information-sufficient level $\rho = 0.80$ [31, 32]. Then the per-emitter SNR in dB is

$$\mathrm{SNR}_m = 10\log_{10}\gamma_m - 20\log_{10}\sigma - 11.02, \tag{28}$$

which is Eq. (11) of Ref. [32] with $\gamma$ replaced by $\gamma_m$.

For the 3D astigmatic PSF with depth-dependent widths $\sigma_x(z)$ and $\sigma_y(z)$ of Appendix C, following Ref. [32], the figure of merit is

$$\eta_m = \frac{\eta_s}{\sigma_x(z_m)\,\sigma_y(z_m)}, \tag{29}$$

then the per-emitter SNR in dB is

$$\mathrm{SNR}_m = 10\log_{10}\gamma_m - 10\log_{10}\left[\sigma_x(z_m)\,\sigma_y(z_m)\right] - 11.02, \tag{30}$$

which is Eq. (9) of Ref. [32] with $\gamma$ replaced by $\gamma_m$.

Under uniform noise and a common intensity, $\bar{n}_m$ reduces to $b + G$, and Eqs. (28) and (30) reduce to the SNR of those two papers. Each frame is reported by the range and mean of $\mathrm{SNR}_m$ over the emitters.

In all simulations, the EM-GML that employs Newton or gradient ascent achieves the same numerical results, so only the results using Newton are reported below.

### *5.3 Asymptotic in multiple frames*

**FOV and camera.** The FOV is $d = 2$ dimensional with size $L_x = L_y = 2400$ nm. The camera has frame size $K_x = K_y = 24$, that is, 576 pixels, with pixel size $\Delta_x = \Delta_y = 100$ nm and exposure time $\Delta_t = 0.01$ s.

**Emitters.** $M = 80$ emitters are placed uniformly at random in the rectangle $[200, 2200] \times [200, 2200]$ of area 4 μm$^2$, i.e., the central block of $20 \times 20$ pixels. The emitters are under a minimum-separation constraint, any candidate within $\delta = 120$ nm of an already accepted position being rejected and redrawn, so the configuration is $\delta$-separated and resolvable while remaining dense with emitter density $80/4 = 20$ emitters/μm$^2$ very high in SMLM. The positions are generated with a fixed Python seed for reproducibility. The emitter intensities $I_m$ are drawn uniformly in $[250000, 350000]$ photons/s.

**PSF.** Given in Appendix C, the *PSF* is the 2D Gaussian of width $\sigma = 108.81$ nm, equal to the approximated standard deviation of an Airy PSF with emitter Alexa 700 and numerical aperture $n_a = 1.4$ [33].

**Noise.** The spatiotemporal density mean of the background Poisson noise $b_{\boldsymbol{k}} \in [4, 6]$ photons/s/nm$^2$ is drawn uniformly. The spatiotemporal density mean and variance of the Gaussian readout noise $\mu_{\boldsymbol{k}} = G_{\boldsymbol{k}} \in [2, 4]$ photons/s/nm$^2$ are drawn uniformly and set equal so that the Poisson approximation of Section 2 applies.

**SNR.** Under these settings the per-emitter information-sufficient SNR is computed by Eq. (28) and reported as the range and mean over the emitters. It is fixed across the frame sweep since summing frames raises the total information without changing the per-frame SNR. The number of frames $N$ is swept from 1 to 1000.

**EM-GML and UGIA-F for $N$ frames.** Because the $N$ frames share the same mean $v(\boldsymbol{k})$, their sum $V_+(\boldsymbol{k}) = \sum_{n=1}^{N} V^{(n)}(\boldsymbol{k}) \sim \text{Poisson}(Nv(\boldsymbol{k}))$ is a sufficient statistic, so the $N$-frame GML is obtained by running the EM algorithm once on $V_+$ rather than on each frame separately. Summing drives the per-pixel count to $Nv(\boldsymbol{k})$, so the noise relative to the signal falls as $1/\sqrt{N}$ and the estimate concentrates at $\boldsymbol{\theta}_0$, the same information-growth mechanism as the high-intensity limit. In the E-step the splitting fraction is unchanged, $\hat{S}_m(\boldsymbol{k}) = (s_m(\boldsymbol{k})/v(\boldsymbol{k}))V_+(\boldsymbol{k})$, and in the M-step the single-emitter likelihood is maximized with mean scaled by $N$,

$$\hat{\boldsymbol{\theta}}_m = \arg\max_{\boldsymbol{\theta}_m} \sum_{\boldsymbol{k}\in\Omega} \left[\hat{S}_m(\boldsymbol{k}) \ln s_m(\boldsymbol{k}) - N s_m(\boldsymbol{k})\right]. \tag{31}$$

The EM-GML uses the small-jitter initialization of Section 5.2. The UGIA-F estimator is also performed on the summed frame. Because $V_+$ carries $N$ times the information of one frame, its FIM is $N\mathbf{F}(\boldsymbol{\theta}_0)$, and the UGIA-F sample is drawn as $\hat{\boldsymbol{\theta}}^F \sim \mathcal{N}\left(\boldsymbol{\theta}_0, \left(N\mathbf{F}(\boldsymbol{\theta}_0)\right)^{-1}\right)$.

**RMSE.** The estimated positions of EM-GML are in known one-to-one correspondence with the true positions through the near-true initialization, and the indices of the UGIA-F sample are known by construction, so no assignment by a partition algorithm [33] is needed. The RMSE for the two estimators is then computed by Eq. (25).

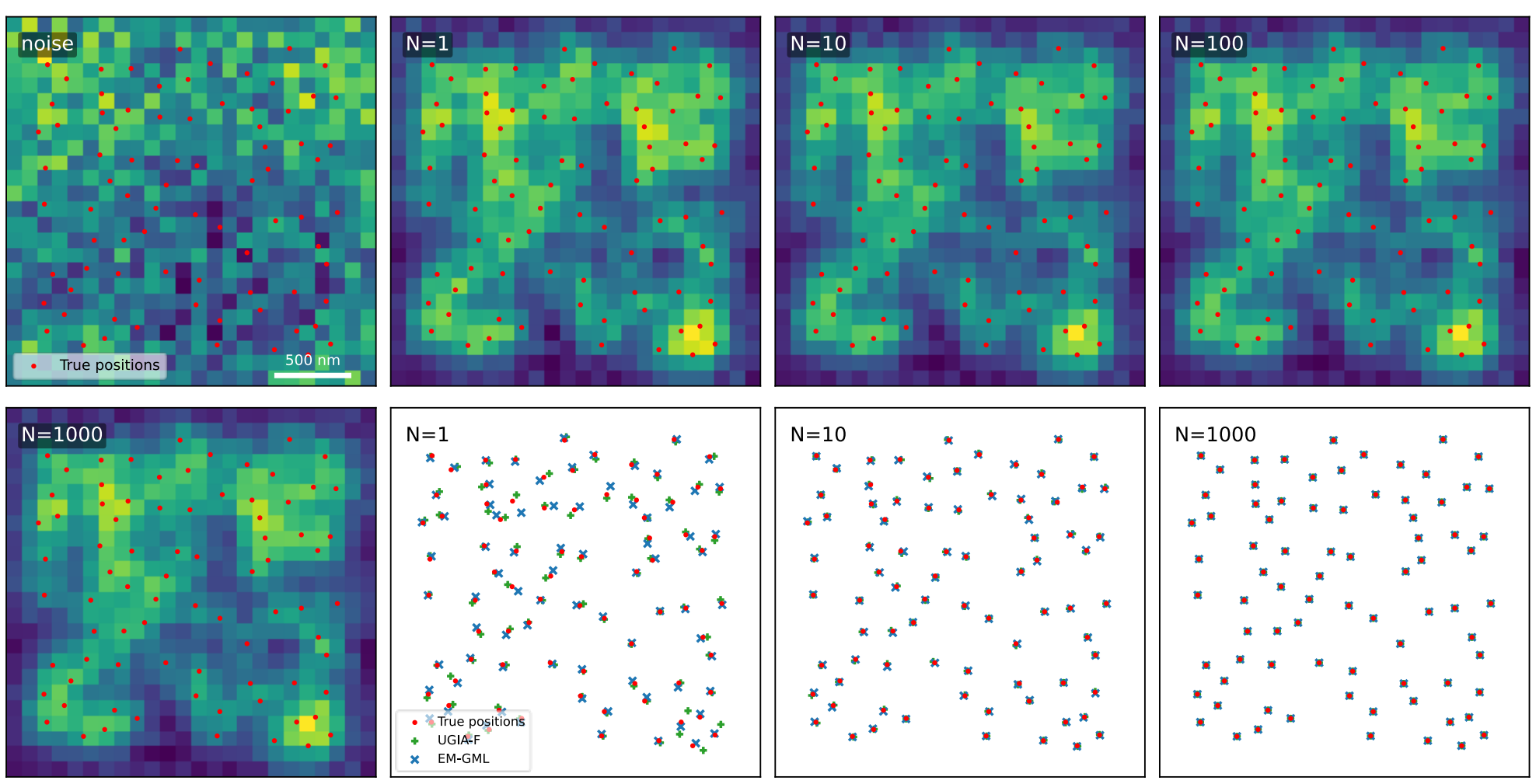

**Fig. 1.** Dense but resolvable 2D configuration of $M = 80$ emitters at density 20 emitters/μm$^2$ and minimum separation $\delta = 120$ nm. Top row: a pure-noise frame, that is, background and readout only, and the summed frames at $N = 1, 10, 100$. Bottom row: the summed frame at $N = 1000$ and the estimated positions of EM-GML and UGIA-F at $N = 1, 10, 1000$. The red dots are the true positions shown in every panel. Each frame is autoscaled, so the noise washes out as $N$ grows while the overlapped pattern stays visually unresolved. The scale bar in the noise frame is 500 nm.

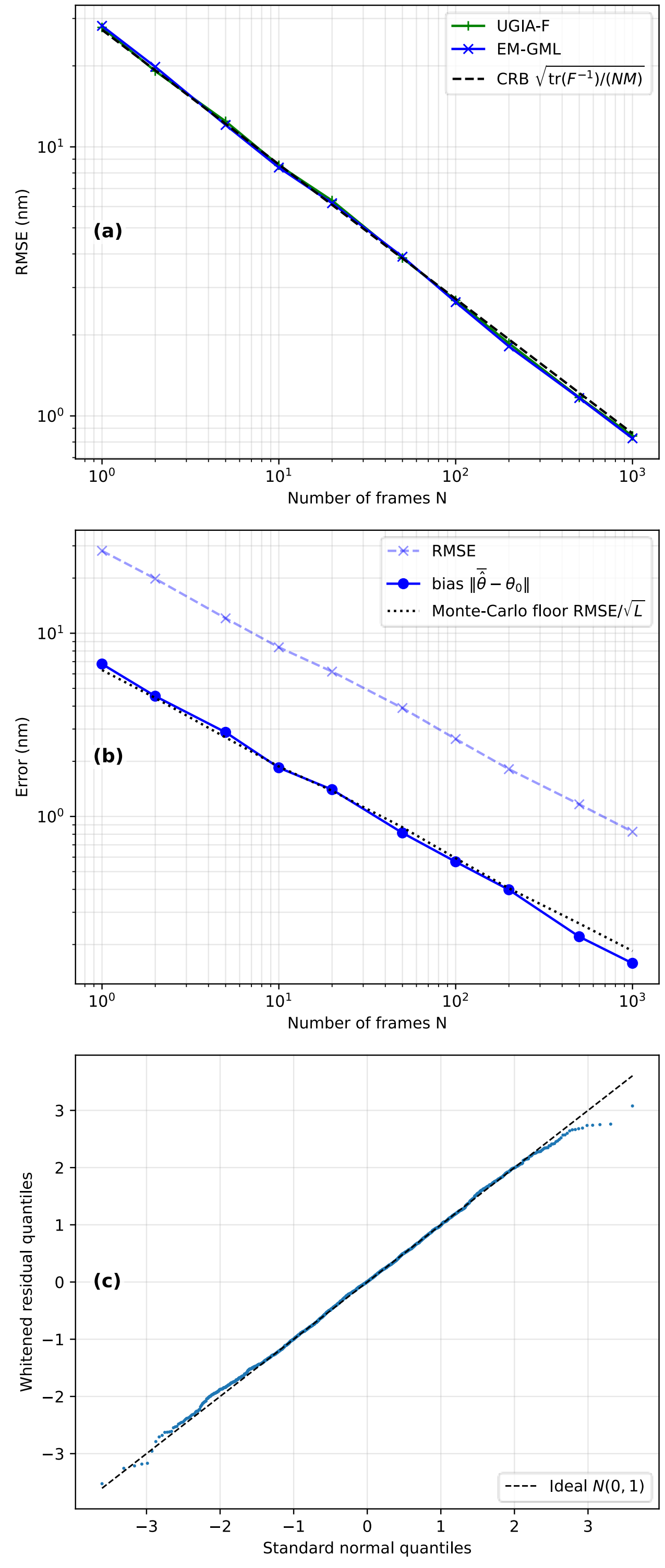

**Fig. 2.** RMSE, bias, and normality of EM-GML versus the number of frames $N$. (a) RMSE versus $N$ on log-log axes for EM-GML and UGIA-F with the $N$-frame CRB $(\mathrm{tr}(\mathbf{F}^{-1})/NM)^{1/2}$

where $\mathbf{F}$ is the one-frame FIM. The three curves coincide and fall as $1/\sqrt{N}$, so EM-GML attains the bound and matches UGIA-F across the sweep. (b) Bias $\left\|\bar{\hat{\boldsymbol{\theta}}} - \boldsymbol{\theta}_0\right\|$ versus $N$ with the EM-GML RMSE of (a) and the Monte-Carlo floor RMSE$/\sqrt{L}$ where $L = 20$ is the number of Monte-Carlo replicates. The measured bias rides the floor, so the true bias is below it and the RMSE is variance-dominated. (c) QQ-plot of the pooled whitened residuals $\mathbf{z} = (N\mathbf{F})^{1/2}(\hat{\boldsymbol{\theta}} - \boldsymbol{\theta}_0)$ against the standard normal at $N = 1000$ where the EM-GML RMSE meets the CRB. The points follow the ideal $\mathcal{N}(0,1)$ line. Moreover, the per-coordinate ratio of the empirical variance to the squared CRLB of Eq. (10) are 0.944 and 0.983 for $x$ and $y$, respectively, both near 1, which confirms asymptotic normality and per-coordinate efficiency.

**Figures.** Fig. 1 shows qualitatively that the summed frame stays an overlapped blur at every $N$ while both estimators converge onto the true positions as $N$ grows, that is, the consistency and the EM-GML and UGIA-F coincidence in position space. Fig. 2 (a) verifies the consistency, the trace efficiency, and the coincidence of EM-GML and UGIA-F in the RMSE as $N$ grows. Fig. 2 (b) verifies the asymptotic unbiasedness, the bias staying at or below the Monte-Carlo floor RMSE$/\sqrt{L}$ and hence negligible relative to the standard deviation. Fig. 2 (c) verifies the asymptotic normality and the per-coordinate efficiency as the variance ratios approach 1, and since UGIA-F is by construction the draw $\mathcal{N}(\boldsymbol{\theta}_0, (N\mathbf{F})^{-1})$, the normality of the whitened residual raises the EM-GML and UGIA-F coincidence from equal RMSE to equal distribution. Together Figs. 1 and 2 verify properties 1-4 of Section 5.1 as $N$ grows.

### *5.4 Asymptotic in emitter intensity*

**FOV and camera.** This simulation verifies Theorem 2 in 3D as the average intensity $I$ grows with a single frame $N = 1$. The FOV is $d = 3$ dimensional with lateral size $L_x = L_y = 2400$ nm and axial size $L_z = 400$ nm or depth $2L_z = 800$ nm. The camera, pixel size, and exposure are as in Section 5.3.

**Emitters.** $M = 40$ emitters are placed randomly in the cuboid $[200, 2200] \times [200, 2200] \times [-400, 400]$ by rejection sampling with 3D minimum separation $\delta = 250$ nm, dense yet resolvable, so the lateral density is 10 emitters/μm$^2$, generated with a fixed Python seed. The intensity weights $\beta_m$ are fixed by drawing $M$ values uniformly in $[250000, 350000]$ photons/s and normalizing to $M^{-1}\sum_{m=1}^{M}\beta_m = 1$, and the average intensity $I$ is then swept so that $I_m = \beta_m I$.

**PSF.** The PSF is the 3D astigmatic Gaussian PSF of Refs. [19, 5] as given in Appendix C and implemented in SMLM_Lib [34], whose lateral widths $\sigma_x(z)$ and $\sigma_y(z)$ vary with depth so that the axial gradient $\partial v(\boldsymbol{k})/\partial z_m$ is nonzero and linearly independent of the lateral gradients, by which C3 holds in $z$. Its calibration parameters are $c = 205$ nm, $a = 290$ nm, $\sigma_{x0} = 140$ nm, $A_x = 0.05$, $B_x = 0.03$, $\sigma_{y0} = 135$ nm, $A_y = -0.01$, and $B_y = 0.02$ as usually used in practice [19, 5].

**Noise.** The background and readout noise are as in Section 5.3. The average intensity $I$ is swept from $10^5$ to $10^{10}$ photons/s, equivalently the per-emitter photon budget $I\Delta_t$ grows, realizing the $I \to \infty$ regime. The per-emitter information-sufficient SNR is computed by Eq. (30) and reported as the range and mean over the emitters, increasing with $I$.

**EM-GML and UGIA-F.** For each $I$, one frame is generated. The EM-GML is run with the small-jitter initialization of Section 5.2, now jittering all three coordinates; and the UGIA-F sample is drawn as $\hat{\boldsymbol{\theta}}^F \sim \mathcal{N}(\boldsymbol{\theta}_0, \mathbf{F}(I)^{-1})$ where $\mathbf{F}(I)$ is the one-frame FIM at intensity $I$ in Eq. (36). The correspondence of estimated and true positions is known for both estimators, so the RMSE is computed by Eq. (25) now over the three coordinates of each position.

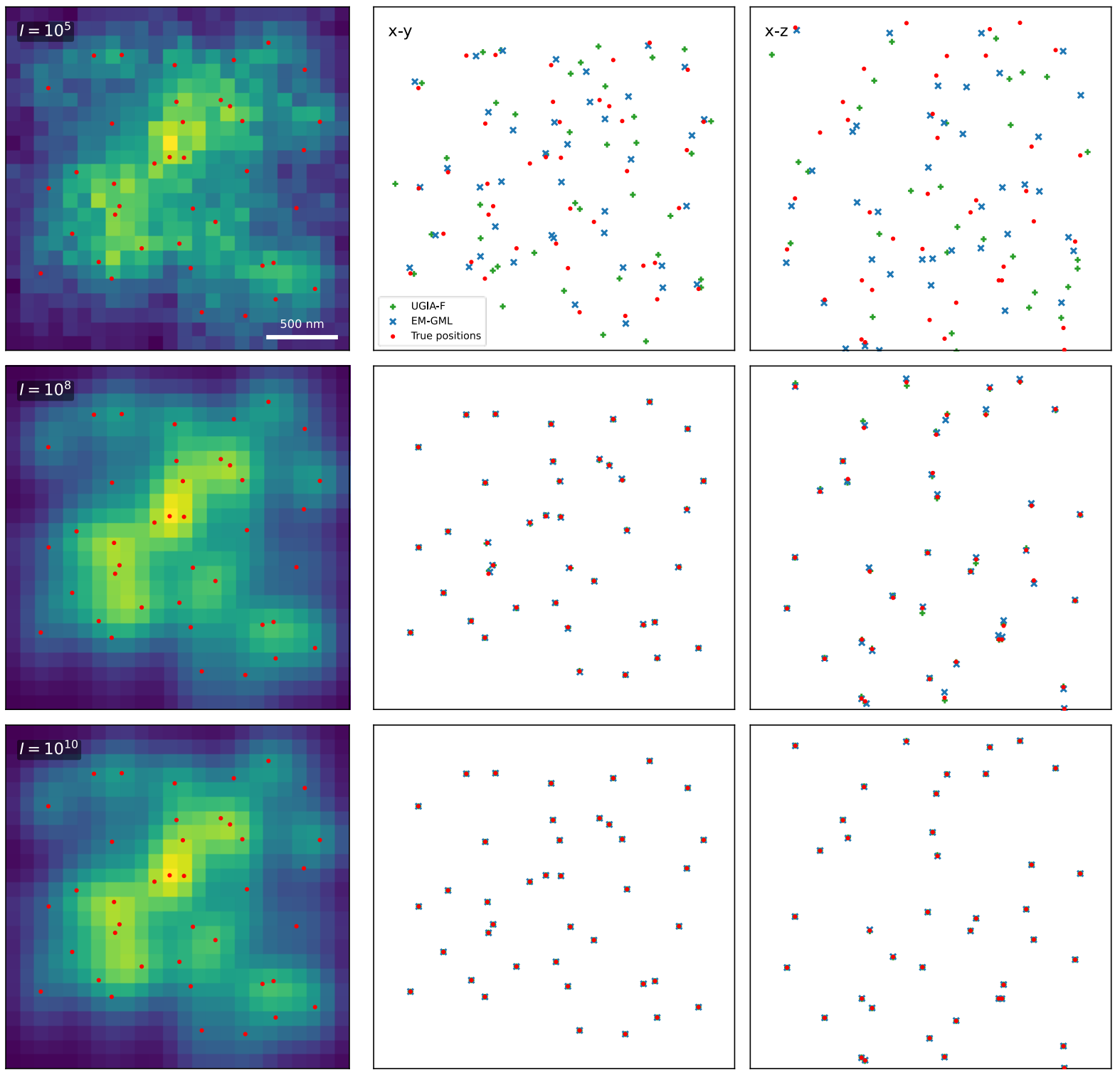

**Fig. 3.** Dense but resolvable 3D astigmatic configuration of $M = 40$ emitters at lateral density 10 emitters/μm² and 3D minimum separation $\delta = 250$ nm. Rows are $I = 10^5, 10^8, 10^{10}$ photons/s. The first column is the data frame, the second the $x$-$y$ reconstruction, and the third the $x$-$z$ reconstruction, with EM-GML and UGIA-F estimates over the true red dots shown in every panel. Each frame is autoscaled, so eventually the noise washes out and the astigmatic spots sharpen as $I$ grows. The scale bar in the $I = 10^5$ frame is 500 nm. The axial range is $[-400, 400]$ nm.

**Figures.** Fig. 3 shows qualitatively that the data frame stays an overlapped astigmatic pattern while both estimators converge onto the true positions in $x$, $y$, and $z$ as $I$ grows, that is, the consistency and the EM-GML and UGIA-F coincidence in 3D position space. Fig. 4 (a) verifies the single-frame trace efficiency and the RMSE coincidence of EM-GML and UGIA-F, lateral and axial alike, EM-GML rising from below the unbiased bound onto the CRB as $I \to \infty$. Fig. 4 (b) verifies the asymptotic unbiasedness, the bias staying at or below the Monte-Carlo floor $\mathrm{RMSE}/\sqrt{L}$ and hence negligible relative to the standard deviation. Fig. 4 (c) verifies the asymptotic normality and the per-coordinate efficiency in 3D, the variance ratios near 0.85 at the highest intensity. The small shortfall is the residual bias of Hypothesis 1 not yet vanished. With the normality the EM-GML and UGIA-F coincidence rises from equal RMSE to equal distribution. Together Figs. 3 and 4 verify properties 1-4 of Section 5.1 in 3D as $I \to \infty$.

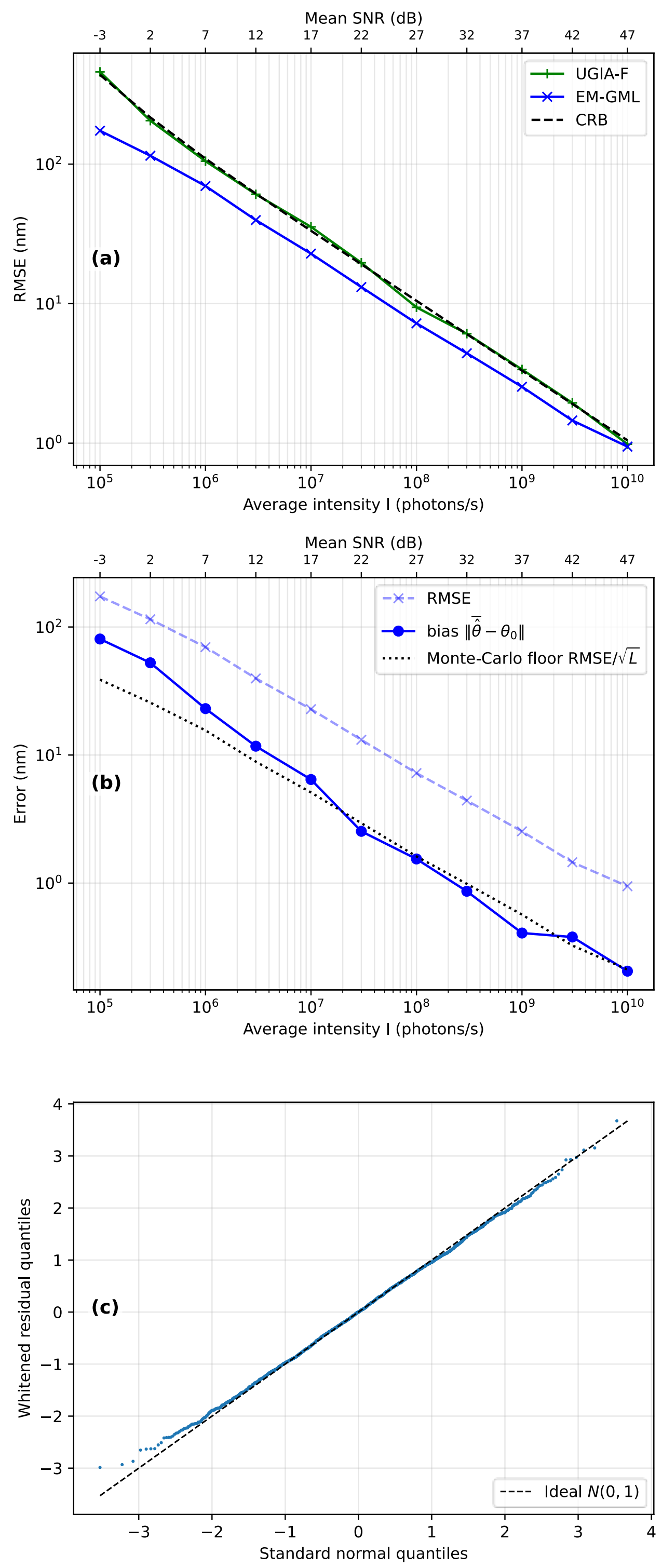

**Fig. 4.** RMSE, bias, and normality of EM-GML versus the average intensity $I$. The top axis of (a) and (b) maps $I$ to the mean SNR. (a) Total RMSE in $x$, $y$, and $z$ for EM-GML and UGIA-F

against the CRB. The lateral RMSE in $x$, $y$ and the axial RMSE in $z$ are nearly equal, so only the total is shown. The EM-GML descends from well below the unbiased bound at low $I$, which is the advantage of Hypothesis 1, and rises onto the CRB and UGIA-F as $I$ grows, which shows the Theorem 2 efficiency. (b) Bias $\left\|\bar{\hat{\boldsymbol{\theta}}} - \boldsymbol{\theta}_0\right\|$ versus $I$ on log-log axes with the EM-GML RMSE of (a) and the Monte-Carlo floor RMSE/$\sqrt{L}$ where $L = 20$ is the number of Monte-Carlo replicates. The measured bias rides the floor, so the true bias is below it and the RMSE is variance-dominated. (c) QQ-plot of the pooled whitened residuals $\mathbf{z} = \mathbf{F}(I)^{1/2}(\hat{\boldsymbol{\theta}} - \boldsymbol{\theta}_0)$ against the standard normal at $I = 10^{10}$ where the EM-GML RMSE meets the CRB. The points follow the ideal $\mathcal{N}(0,1)$ line, and the per-coordinate variance ratios are 0.888, 0.792 and 0.857 for $x$, $y$, and $z$, respectively, all near 0.85, which confirms asymptotic normality and the approach to per-coordinate efficiency.

### *5.5 EM-GML outperformance over CRB in the finite-photon regime*

**Setup.** The configuration is that of Section 5.3 with a single frame $N = 1$ and the 2D Gaussian PSF in Appendix C. The number of emitters $M$ is swept from 1 to 200 inside the fixed central block of Section 5.3, so the emitter density rises up to $200/4 = 50$ emitters/μm$^2$, very high in SMLM. For each $M$ the positions are $\delta$-separated by rejection sampling with $\delta = 100$ nm and a fixed Python seed, the intensities are drawn uniformly in $[250000, 350000]$ photons/s, and the background and readout noises are as in Section 5.3.

**EM-GML and UGIA-F.** For each $M$ one frame is generated. The EM-GML is run with the small-jitter initialization of Section 5.2 and the UGIA-F sample is drawn as $\hat{\boldsymbol{\theta}}^F \sim \mathcal{N}(\boldsymbol{\theta}_0, \mathbf{F}(\boldsymbol{\theta}_0)^{-1})$. The correspondence between the estimated and true positions is known for both estimators, so the RMSE is computed by Eq. (25).

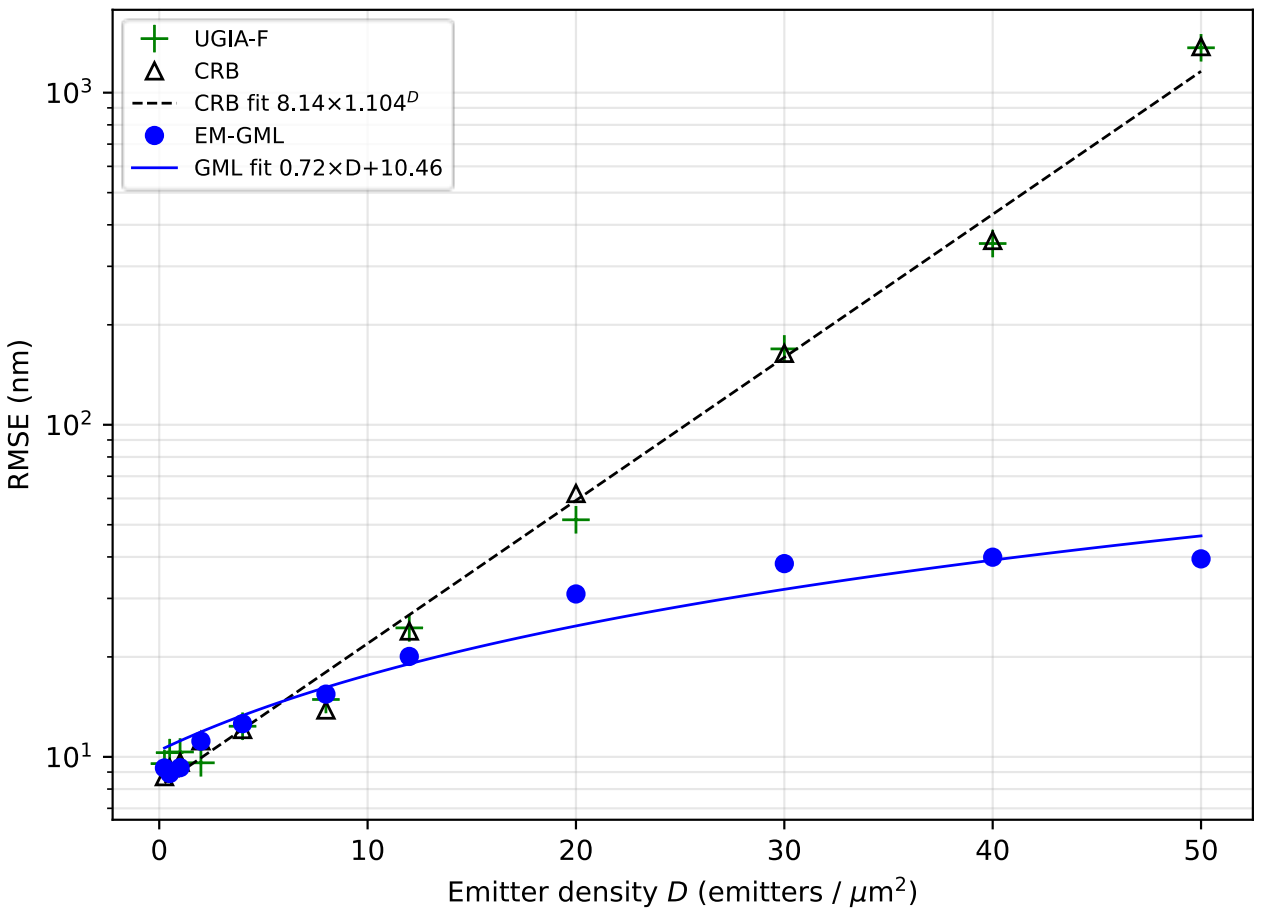

**Fig. 5.** On a semi-logarithmic axis, computed with $L = 20$ Monte-Carlo replicates, RMSE versus emitter density $D$ for EM-GML and UGIA-F with the frame CRB reference $(\mathrm{tr}(\mathbf{F}^{-1})/M)^{1/2}$ at a single frame $N = 1$. The emitter density $D$ is swept over $M \in \{1, 4, 8, 16, 32, 48, 80, 120, 160, 200\}$ inside the fixed central block. EM-GML is the filled blue circles, UGIA-F the green plus signs, and the single-frame CRB the open black triangles. As $D$ grows the CRB and the UGIA-F benchmark increase exponentially, well fit by $8.14 \times 1.104^D$ nm, while the EM-GML RMSE increases only roughly linearly, well fit by $0.72\,D + 10.46$ nm. The biased EM-GML therefore falls increasingly far below the unbiased benchmark as the density rises.

**Figures.** As shown in Fig. 5, the linear-versus-exponential separation is the gap of Hypothesis 1 in Section 3.6, here exposed by raising the density at a single frame rather than by lowering the photon count. It reproduces the finite-photon behavior of Ref. [24] where the biased EM-GML estimator outperforms the unbiased CRB once the emitters are dense enough that the Fisher information becomes ill-conditioned. Because UGIA-F tracks the CRB

throughout, the gap is not an artifact of the benchmark but the genuine advantage of the biased GML in the finite-information regime. The substantial gap between the EM-GML and UIGA-F/CRB in the high-density regime implies a great potential to develop advanced localization algorithms for high-throughput super spatiotemporal resolution SMLM.

## 6. Conclusions

This paper established the asymptotic efficiency of the GML estimator for the multi-emitter SMLM imaging model. For any fixed configuration of resolvable emitters in 2D or 3D, the GML estimator of the emitter positions is consistent, asymptotically normal, asymptotically unbiased, and asymptotically efficient, attaining the Cramér–Rao bound as the number of frames grows in Theorem 1 and as the average emitter intensity grows in Theorem 2. Corollary 1 states that in either limit the GML estimator and the unbiased UGIA-F benchmark converge to the same distribution, so the two coincide wherever the information is large.

The asymptotic machinery is classical, so the contribution is its precise application to the SMLM Poisson model. This required defining the GML estimator for the model, verifying the regularity conditions, in particular the strict positivity of the Fisher information that keeps the model cleanly regular for resolvable emitters, and then invoking the classical theory. The theory was verified by simulation with the expectation-maximization realization EM-GML, Theorem 1 with a 2D Gaussian point spread function and Theorem 2 with a 3D astigmatic one, in both of which EM-GML and UGIA-F approach the Cramér–Rao bound together.

The results are restricted to fixed, resolvable configurations. As the separation $\delta$ tends to zero the Fisher information becomes singular, its inverse diverges, and asymptotic efficiency ceases to be informative, which is the high-density regime where a biased estimator can do better. There, in the finite-photon regime, the biased EM-GML attains an error below the unbiased Cramér–Rao bound and below UGIA-F, a gap that widens with emitter density, the high-density advantage studied as Hypothesis 1. This gap has two consequences for practice. First, EM-GML together with UGIA-F provides a benchmark for localization algorithms that do not know the system parameters, including neural networks. Second, the gap motivates the design of advanced localization algorithms, including neural networks, that exploit it at high emitter density to achieve high-throughput super spatiotemporal resolution SMLM.

## Appendix A: Proof of Theorems 1 and 2

### *A.1 Theorem 1 (i): Consistency*

The frames are i.i.d. with continuous, identifiable likelihood on the compact set $\Theta$. Define the Kullback–Leibler contrast

$$\mathbb{K}(\boldsymbol{\theta}_0, \boldsymbol{\theta}) = \mathbb{E}_{\boldsymbol{\theta}_0}[\ell_1(\boldsymbol{\theta}_0) - \ell_1(\boldsymbol{\theta})] = \sum_{\boldsymbol{k}\in\Omega}\left[v(\boldsymbol{k}) - v_0(\boldsymbol{k}) - v_0(\boldsymbol{k})\ln\frac{v(\boldsymbol{k})}{v_0(\boldsymbol{k})}\right] \geq 0, \quad (32)$$

where $v_0(\boldsymbol{k}) = v(\boldsymbol{k})|_{\boldsymbol{\theta}_0}$. Each summand is the Kullback-Leibler divergence between the Poisson laws of means $v_0(\boldsymbol{k})$ and $v(\boldsymbol{k})$, so $\mathbb{K}(\boldsymbol{\theta}_0, \boldsymbol{\theta})$ is nonnegative and continuous (C2), and zero only at $\boldsymbol{\theta} = \boldsymbol{\theta}_0$ up to emitter relabeling since $v(\boldsymbol{k}) \geq v_{\min} > 0$; the local identifiability from C3 singles out $\boldsymbol{\theta}_0$ in a neighborhood. By the uniform law of large numbers, from compactness C1 and domination from C2,

$$\sup_{\boldsymbol{\theta}\in\Theta}\left|\frac{1}{N}\ell_N(\boldsymbol{\theta}) - \mathbb{E}_{\boldsymbol{\theta}_0}[\ell_1(\boldsymbol{\theta})]\right| \xrightarrow{a.s.} 0. \quad (33)$$

The maximizer of a uniformly convergent sequence of objectives converges to the maximizer of the limit, i.e., $\hat{\boldsymbol{\theta}}_N^{\mathrm{GML}} \to \boldsymbol{\theta}_0$ [12]. Because the global maximizer is used, spurious *local* maxima, which have strictly smaller expected log-likelihood, are excluded in the limit.

### *A.2 Theorem 1 (ii): Asymptotic normality and efficiency*

Since $\boldsymbol{\theta}_0 \in \operatorname{int}\Theta$ and $\hat{\boldsymbol{\theta}}_N^{\mathrm{GML}} \to \boldsymbol{\theta}_0$, eventually the total score $\boldsymbol{s}_N(\hat{\boldsymbol{\theta}}_N^{\mathrm{GML}}) = \mathbf{0}$. A first-order Taylor expansion of the score of Eq. (8) about $\boldsymbol{\theta}_0$ gives for some $\overline{\boldsymbol{\theta}}$ on the segment to $\boldsymbol{\theta}_0$,

$$\frac{1}{\sqrt{N}}\boldsymbol{s}_N(\boldsymbol{\theta}_0) + \left[\frac{1}{N}\nabla_{\boldsymbol{\theta}}\boldsymbol{s}_N(\overline{\boldsymbol{\theta}})\right]\sqrt{N}\left(\hat{\boldsymbol{\theta}}_N^{\mathrm{GML}} - \boldsymbol{\theta}_0\right) = \mathbf{0}. \tag{34}$$

We have two limits.

- *Score central limit theorem.* $\boldsymbol{s}_1(\boldsymbol{\theta}_0)$ has mean $\mathbf{0}$ and covariance $\mathbf{F}(\boldsymbol{\theta}_0)$, so by the central limit theorem $\boldsymbol{s}_N(\boldsymbol{\theta}_0)/\sqrt{N} \xrightarrow{d} \mathcal{N}(\mathbf{0}, \mathbf{F}(\boldsymbol{\theta}_0))$.
- *Observed information.* By the law of large numbers and C4, which controls the remainder via the third derivative, $-N^{-1}\nabla_{\boldsymbol{\theta}}\boldsymbol{s}_N(\overline{\boldsymbol{\theta}}) \xrightarrow{p} \mathbf{F}(\boldsymbol{\theta}_0)$.

Slutsky's theorem yields

$$\sqrt{N}\left(\hat{\boldsymbol{\theta}}_N^{\mathrm{GML}} - \boldsymbol{\theta}_0\right) \xrightarrow{d} \mathbf{F}(\boldsymbol{\theta}_0)^{-1}\,\mathcal{N}(\mathbf{0}, \mathbf{F}(\boldsymbol{\theta}_0)) = \mathcal{N}(\mathbf{0}, \mathbf{F}(\boldsymbol{\theta}_0)^{-1}). \tag{35}$$

The asymptotic covariance equals the CRB $\mathbf{F}^{-1}$, which is the definition of asymptotic *efficiency* [14, 17, 16].

### A.3 Theorem 1 (iii): Asymptotic unbiasedness

$\hat{\boldsymbol{\theta}}_N^{\mathrm{GML}}$ lies in the compact set $\Theta$, hence is bounded and uniformly integrable; the convergence in probability of A.1 therefore upgrades to convergence in mean, $\mathbb{E}\left(\left\|\hat{\boldsymbol{\theta}}_N^{\mathrm{GML}} - \boldsymbol{\theta}_0\right\|\right) \to 0$. The standard higher-order stochastic expansion of the MLE [17] gives bias $= O(N^{-1})$, dominated by the $O(N^{-1/2})$ standard deviation; the estimator is thus asymptotically unbiased and efficient.

### A.4 Theorem 2: $I \to \infty$ regime

Let $I \to \infty$ with the weights $\beta_m$ fixed and so $I_m = I\beta_m$. Then $v(\boldsymbol{k}) = \Delta_t\Delta_x\Delta_y\, I\, Q(\boldsymbol{k}) + b(\boldsymbol{k})$, the signal grows with $I$ while the mixed-noise mean $b(\boldsymbol{k})$ stays fixed, and the FIM is an explicit function of $I$ from Eq. (9),

$$\mathbf{F}(I) = \sum_{\boldsymbol{k}\in\Omega} \frac{\left(\Delta_t\Delta_x\Delta_y\, I\right)^2}{\Delta_t\Delta_x\Delta_y\, I\, Q(\boldsymbol{k}) + b(\boldsymbol{k})} \frac{\partial Q(\boldsymbol{k})}{\partial\boldsymbol{\theta}} \frac{\partial Q(\boldsymbol{k})}{\partial\boldsymbol{\theta}}^{\top} \xrightarrow{I\to\infty} I\,\mathbf{C} \tag{36}$$

where

$$\mathbf{C} = \Delta_t\Delta_x\Delta_y \sum_{\boldsymbol{k}:\, Q(\boldsymbol{k})>0} \frac{1}{Q(\boldsymbol{k})} \frac{\partial Q(\boldsymbol{k})}{\partial\boldsymbol{\theta}} \frac{\partial Q(\boldsymbol{k})}{\partial\boldsymbol{\theta}}^{\top} \succ \mathbf{0} \tag{37}$$

because the fixed background $b(\boldsymbol{k})$ is negligible beside the growing signal. Then every eigenvalue grows like $I$ and in particular $\lambda_{\min}(\mathbf{F}(I)) \to \infty$, equivalently $\mathbf{F}(I)^{-1} \to \mathbf{0}$; the CRB shrinks as $1/I$, recovering the familiar photon-limited $R \propto 1/\sqrt{I}$ scaling. Write $F_V^I(\boldsymbol{\theta})$ for the law (i.e., probability distribution) and $f_V^I(\boldsymbol{\theta})$ for the density of the frame $V$ under parameter $\boldsymbol{\theta}$ at average intensity $I$. The single-frame family $\{F_V^I(\boldsymbol{\theta}) : \boldsymbol{\theta} \in \Theta\}$ is then locally asymptotically normal at $\boldsymbol{\theta}_0$ with norming matrix $\mathbf{F}(I)^{1/2}$. For a local parameter $\boldsymbol{h} \in \mathbb{R}^{dM}$, the log-likelihood ratio of measures in terms of the Radon-Nikodym derivative is

$$\ln\frac{dF_V^I\left(\boldsymbol{\theta}_0 + \mathbf{F}(I)^{-1/2}\boldsymbol{h}\right)}{dF_V^I(\boldsymbol{\theta}_0)} = \ln\frac{f_V^I\left(\boldsymbol{\theta}_0 + \mathbf{F}(I)^{-1/2}\boldsymbol{h}\right)}{f_V^I(\boldsymbol{\theta}_0)} = \boldsymbol{h}^{\top}\underbrace{\mathbf{F}(I)^{-1/2}\dot{\ell}_I(\boldsymbol{\theta}_0)}_{\boldsymbol{Z}_I}$$

$$+\frac{1}{2}\boldsymbol{h}^{\top}\left[\mathbf{F}(I)^{-1/2}\ddot{\ell}_I(\boldsymbol{\theta}_0)\,\mathbf{F}(I)^{-1/2}\right]\boldsymbol{h} + R_I = \boldsymbol{h}^{\top}\boldsymbol{Z}_I - \frac{1}{2}\|\boldsymbol{h}\|^2 + o_p(1). \tag{38}$$

Here $\dot{\ell}_I(\boldsymbol{\theta}_0)$ is the score $\boldsymbol{s}_1(\boldsymbol{\theta}_0)$ of Eq. (8), $\ddot{\ell}_I(\boldsymbol{\theta}_0)$ its Hessian, and $R_I$ the third-order Taylor remainder, all at average intensity $I$. The normalized score $\boldsymbol{Z}_I \xrightarrow{d} \mathcal{N}(\mathbf{0}, \mathbf{I}_{dM})$ because the Poisson log-likelihood is smooth with mean bounded below (C2) and the third-order remainder is controlled by $\mathbf{F}(I)^{-1/2}$ (C4). Le Cam's theory [15, 17] then gives efficiency of the GML, namely $\mathbf{F}(I)^{1/2}\left(\hat{\boldsymbol{\theta}}^{\mathrm{GML}} - \boldsymbol{\theta}_0\right) \xrightarrow{d} \mathcal{N}(\mathbf{0}, \mathbf{I}_{dM})$.

## Appendix B: Noise-free single-emitter ML estimator

In the M-step of Eq. (24), the position of emitter $m$ is estimated from the noise-free estimated image $\hat{S}_m^{(l+1)}(\boldsymbol{k})$. Two iterative solvers are used, gradient ascent and Newton. Both run on an inner index $n$ from the initial $\boldsymbol{\theta}_m(0) = \boldsymbol{\theta}_m^{(l)}$ and return $\hat{\boldsymbol{\theta}}_m^{(l+1)}$.

**Gradient ascent.** The gradient $\boldsymbol{L}_m^{(l+1)}(\boldsymbol{\theta}_m)$ of the noise-free single-emitter log-likelihood $\ln f_{S_m}\left(\hat{S}_m^{(l+1)}; \boldsymbol{\theta}_m\right)$ has the components for $j = 1, \dots, d$,

$$\left[\boldsymbol{L}_m^{(l+1)}(\boldsymbol{\theta}_m)\right]_j = \frac{\partial}{\partial \theta_{mj}} \ln f_{S_m}\left(\hat{S}_m^{(l+1)}; \boldsymbol{\theta}_m\right) = \Delta_t \Delta_x \Delta_y \, I_m \sum_{\boldsymbol{k}\in\Omega} \left(\frac{\hat{S}_m^{(l+1)}(\boldsymbol{k})}{s_m(\boldsymbol{k})} - 1\right) \frac{\partial q_m(\boldsymbol{k})}{\partial \theta_{mj}}. \quad (39)$$

The position is updated along the gradient,

$$\boldsymbol{\theta}_m(n+1) = \boldsymbol{\theta}_m(n) + \alpha\, \boldsymbol{L}_m^{(l+1)}\big(\boldsymbol{\theta}_m(n)\big), \quad (40)$$

with a constant step $\alpha > 0$, terminating when $\|\boldsymbol{\theta}_m(n+1) - \boldsymbol{\theta}_m(n)\| < \varepsilon$. With the iteration count fixed, all $M$ emitters are updated in parallel as a GPU batch, which speeds the compute by about a factor of $M$ over a serial sweep.

**Newton.** The Newton algorithm uses the Hessian $\mathbf{H}(\boldsymbol{\theta}_m)$ of the noise-free single-emitter log-likelihood with elements for $i, j = 1, \dots, d$,

$$[\mathbf{H}(\boldsymbol{\theta}_m)]_{ij} = \frac{\partial^2}{\partial \theta_{mi}\, \partial \theta_{mj}} \ln f_{S_m}\left(\hat{S}_m^{(l+1)}; \boldsymbol{\theta}_m\right), \quad (41)$$

and updates

$$\boldsymbol{\theta}_m(n+1) = \boldsymbol{\theta}_m(n) - \mathbf{H}^{-1}\big(\boldsymbol{\theta}_m(n)\big)\, \boldsymbol{L}_m^{(l+1)}\big(\boldsymbol{\theta}_m(n)\big). \quad (42)$$

The Hessian in Eq. (41) can be replaced by its expectation, $E[\mathbf{H}(\boldsymbol{\theta}_m)] = -\mathbf{F}(\boldsymbol{\theta}_m)$, the $d \times d$ noise-free single-emitter FIM with elements for $i, j = 1, \dots, d$,

$$[\mathbf{F}(\boldsymbol{\theta}_m)]_{ij} = \Delta_t \Delta_x \Delta_y \, I_m \sum_{\boldsymbol{k}\in\Omega} \frac{1}{q_m(\boldsymbol{k})} \frac{\partial q_m(\boldsymbol{k})}{\partial \theta_{mi}} \frac{\partial q_m(\boldsymbol{k})}{\partial \theta_{mj}}, \quad (43)$$

evaluated at the current estimate. The resulting update is

$$\boldsymbol{\theta}_m(n+1) = \boldsymbol{\theta}_m(n) + \alpha\, \mathbf{F}^{-1}\big(\boldsymbol{\theta}_m(n)\big)\, \boldsymbol{L}_m^{(l+1)}\big(\boldsymbol{\theta}_m(n)\big) \quad (44)$$

with $\alpha = 1$. Near an LML point the Newton algorithm converges quadratically in the estimation error whereas the gradient ascent converges geometrically with a rate limited by the spread of the eigenvalues of the Hessian, which constrains the step $\alpha$.

## Appendix C: The 2D Gaussian and 3D astigmatic PSFs

The 2D simulations of Sections 5.3 and 5.5 use the 2D Gaussian PSF, and the 3D simulation of Section 5.4 uses the 3D astigmatic Gaussian PSF. Both are given here.

**2D Gaussian PSF.** For a 2D Gaussian PSF of width $\sigma$ [18, 5], the PSF density at the camera-plane point $\boldsymbol{u} = (x, y)$ is

$$q_m(\boldsymbol{u}) = \frac{1}{2\pi\sigma^2} \exp\left(-\frac{(x - x_m)^2 + (y - y_m)^2}{2\sigma^2}\right), \quad (45)$$

whose pixel average by Eq. (3) is

$$q_m(\boldsymbol{k}) = \frac{1}{\Delta_x \Delta_y}\left[\Phi\left(\frac{\Delta_x(k_x + 1) - x_m}{\sigma}\right) - \Phi\left(\frac{\Delta_x k_x - x_m}{\sigma}\right)\right]$$
$$\times\left[\Phi\left(\frac{\Delta_y(k_y + 1) - y_m}{\sigma}\right) - \Phi\left(\frac{\Delta_y k_y - y_m}{\sigma}\right)\right], \quad (46)$$

with $\Phi$ the standard normal CDF.

**3D astigmatic PSF.** The 3D astigmatic Gaussian PSF of Refs. [19, 5] has the density for $\boldsymbol{u} = (x, y)$

$$q_m(\boldsymbol{u}) = \frac{1}{2\pi\sigma_x(z_m)\sigma_y(z_m)} \exp\left(-\frac{(x - x_m)^2}{2\sigma_x^2(z_m)} - \frac{(y - y_m)^2}{2\sigma_y^2(z_m)}\right), \quad (47)$$

where $\sigma_x^2(z) = \sigma_{x0}^2(1 + (z + c)^2/a^2 + A_x\,(z + c)^3/a^3 + B_x\,(z + c)^4/a^4)$ and $\sigma_y^2(z) = \sigma_{y0}^2\left(1 + (z - c)^2/a^2 + A_y\,(z - c)^3/a^3 + B_y\,(z - c)^4/a^4\right)$ with focal offset $c$, depth parameter $a$, in-focus widths $\sigma_{x0}$ and $\sigma_{y0}$, and shape coefficients $A_x$, $B_x$, $A_y$, $B_y$. The astigmatism makes $\sigma_x$ and $\sigma_y$ vary oppositely with depth, so the frame encodes the depth $z_m$. The pixel average by Eq. (3) replaces $\sigma$ in Eq. (46) by $\sigma_x(z_m)$ along $x$ and $\sigma_y(z_m)$ along $y$.

### Funding

### Acknowledgements

### Disclosures

The author declares no conflicts of interest.

### Data Availability

Data underlying the results presented in this paper are available in Ref. [35].

### References

1. E. Betzig, G. H. Patterson, R. Sougrat, O. W. Lindwasser, S. Olenych, J. S. Bonifacino, M. W. Davidson, J. Lippincott-Schwartz, and H. F. Hess, "Imaging intracellular fluorescent proteins at nanometer resolution," Science, **313**(5793), 1642–1645(2006).
2. M. J. Rust, M. Bates, and X. Zhuang, "Sub-diffraction-limit imaging by stochastic optical reconstruction microscopy (STORM)," Nat. Methods, **3**(10), 793-796(2006).
3. S. T. Hess, T. P. Girirajan, and M. D. Mason, "Ultra-high resolution imaging by fluorescence photoactivation localization microscopy," Biophys. J., **91**(11), 4258–4272(2006).
4. M. Heilemann , S. v. Linde , M. Schüttpelz , R. Kasper, B. Seefeldt, A. Mukherjee, P. Tinnefeld, and M. Sauer, "Subdiffraction-resolution fluorescence imaging with conventional fluorescent probes," Angew. Chem. Int. Ed., **47**(33), 6172–6176(2008).
5. Y. Sun, "Localization precision of stochastic optical localization nanoscopy using single frames," J. Biomed. Optics, **18**(11), 111418-14(2013).
6. Y. Sun, and Y. Guan, "Effect of unknown emitter intensities on localization accuracy in stochastic optical localization nanoscopy using single frames," JOSA A, **38**(12), 1830-1840(2021).
7. Y. Sun, "Root mean square minimum distance as a quality metric for stochastic optical localization nanoscopy images," Sci. Reports, **8**(1), 17211(2018).
8. R. J. Ober, S. Ram and E. S. Ward, "Localization accuracy in single-molecule microscopy," Biophys. J., **86**(2), 1185-1200(2004).
9. I. K. Mortensen, L. S. Churchman, J. A. Spudich and H. Flyvbjerg, "Optimized localization analysis for single-molecule tracking and super-resolution microscopy," Nat. Methods, **7**(5), 377-381(2010).

10. B. Rieger and S. Stallinga, “The Lateral and Axial Localization Uncertainty in Super-Resolution Light Microscopy,” ChemPhysChem, **15**(4), 664-670(2014).
11. J. Chao, E. S. Ward and R. J. Ober, “Fisher information theory for parameter estimation in single molecule microscopy: tutorial,” JOSA A, **33**(7), B36-B57(2016).
12. A. Wald, “Note on the consistency of the maximum likelihood estimate,” Ann. Math. Statist., **20**(4), 595-601(1949).
13. C. R. Rao, “Information and the accuracy attainable in the estimation of statistical parameters,” Bull. Calcutta Math. Soc., **37**(3), 81-91(1945).
14. H. Cramér, Mathematical methods of statistics. Vol. 9, Princeton, NJ: Princeton University Press (1999).
15. L. Cam, “Locally asymptotically normal families of distributions,” Univ. California Publ. Statist., **3**, 37-98(1960).
16. E. L. Lehmann and G. Casella, Theory of point estimation, New York, NY: Springer (1998).
17. A. W. Van der Vaart, Asymptotic statistics, Cambridge, UK: Cambridge University Press (2000).
18. S. Stallinga and B. Rieger, “Accuracy of the Gaussian point spread function model in 2D localization microscopy,” Opt. Express, **18**(24), 24461-24476(2010).
19. B. Huang, W. Wang, M. Bates and X. Zhuang, “Three-dimensional super-resolution imaging by stochastic optical reconstruction microscopy,” Science, **319**(5864), 810-813(2008).
20. S. R. P. Pavani, M. A. Thompson, J. S. Biteen, S. J. Lord, N. Liu, R. J. Twieg, R. Piestun and W. E. Moerner, “Three-dimensional, single-molecule fluorescence imaging beyond the diffraction limit by using a double-helix point spread function,” Proc. Natl. Acad. Sci. USA, **106**(9), 2995-2999(2009).
21. Y. Shechtman, S. J. Sahl and A. S. Backer, “Optimal point spread function design for 3D imaging,” Phys. Rev. Lett., **113**(13), 133902(2014).
22. Y. Li, M. Mund, P. Hoess, J. Deschamps, U. Matti, B. Nijmeijer, V. J. Sabinina, J. Ellenberg, I. Schoen and J. Ries, “Real-time 3D single-molecule localization using experimental point spread functions," Nat. Methods, **15**(5), 367-369(2018).
23. A. P. Dempster, N. M. Laird and D. B. Rubin, “Maximum likelihood from incomplete data via the EM algorithm,” J. R. Stat. Soc. Ser. B, **39**(1), 1-22(1977).
24. Y. Sun, Y. Zeng, W. Yen, J. M. Tarbell and B. M. Fu, “Expectation maximization can outperform Cramer-Rao lower bound in stochastic optical localization nanoscopy,” Quantitative Bioimaging Conf. (QBI2014) (2014).
25. F. Huang, T. M. Hartwich, F. E. Rivera-Molina, Y. Lin, W. C. Duim, J. J. Long, P. D. Uchil, J. R. Myers, M. A. Baird, W. Mothes, M. W. Davidson, D. Toomre and J. Bewersdorf, “Video-rate nanoscopy using sCMOS camera–specific single-molecule localization algorithms,” Nat. Methods, **10**(7), 653-658(2013).
26. M. Born and E. Wolf, Principles of Optics, 7th Ed ed., Cambridge, UK: Cambridge University Press (1999).
27. S. F. Gibson and F. Lanni, “Experimental test of an analytical model of aberration in an oil-immersion objective lens used in three-dimensional light microscopy,” J. Opt. Soc. Am. A, **8**(10), 1601-1613(1999).
28. B. Richards and E. Wolf, “Electromagnetic diffraction in optical systems, II. Structure of the image field in an aplanatic system,” Proc. R. Soc. A, **253**(1274), 358-379(1959).
29. H. White, “Maximum likelihood estimation of misspecified models,” Econometrica, 1-25(1982).
30. S. M. Ross, Stochastic Processes, 2nd Ed. ed., New York: John Wiley and Sons (1996).
31. Y. Sun, “Information sufficient segmentation and signal-to-noise ratio in stochastic optical localization nanoscopy,” Opt. Lett., **45**(21), 6102-6105(2020).
32. M. Sun and Y. Sun, “Information sufficient segmentation and signal-to-noise ratio for 3D astigmatism stochastic optical localization nanoscopy,” Electron. Lett., **58**(2), 58-60(2022).
33. Y. Sun, “Partition of estimated locations: an approach to accurate quality metrics for stochastic optical localization Nanoscopy,” J. Opt. Soc. Am. A, **39**(12), 2307-2315(2022).
34. Y. Sun, SMLM_Lib. [Online] Available: https://github.com/SunCCNY/Nanoscopy/tree/main/SMLM_Lib. [Accessed July 5, 2026].
35. Y. Sun, GML-EM. [Online] Available: https://github.com/SunCCNY/Nanoscopy/tree/main/2026-JOSAA-GML-EM. [Accessed July 16, 2026].